\documentclass[aos,MSNbibl,dvips]{arximspdf}

\doi{10.1214/12-AOS1043} %kopijuoti is PTS
\volume{40}
\issue{5}
\pubyear{2012}
\firstpage{2601}
\lastpage{2633}

\makeatletter

\newcommand{\rright}{\right}
\newcommand{\lleft}{\left}

\newtheorem{theorem}{Theorem}[section]

\newtheorem{lemma}{Lemma}[section]

\newproclaim{definition}{Definition}[section]

\newcommand{\reals}{\mathbb{R}}
\newcommand{\integers}{\mathbb{Z}}

\newcommand{\pr}{\mathbb{P}}
\newcommand{\ex}{\mathbb{E}}
\newcommand{\normal}{N}

\newcommand{\Obs}{T_{[0,1]^d}}

\newcommand{\convd}{\stackrel{d}{\longrightarrow}}
\newcommand{\convas}{\stackrel{\mathrm{a.s.}}{\longrightarrow}}
\newcommand{\X}{Y}
\newcommand{\Xvec}{\underline{Y}}

\newcommand{\Dvec}{A}
\newcommand{\DDvec}{B}
\newcommand{\epsvec}{\underline{\varepsilon}}
\newcommand{\Evec}{\underline{\xi}}

\newcommand{\um}{C_m}
\newcommand{\lsigma}{c_{\sigma}}
\newcommand{\usigma}{C_{\sigma}}

\newcommand{\sset}{S}

\makeatother

\begin{document}
\begin{frontmatter}

\title{Nonparametric regression for locally stationary time
series\thanksref{T1}}
\runtitle{Locally stationary nonparametric regression}

\thankstext{T1}{Supported by the DFG-SNF research group FOR916.}

\begin{aug}
\author[A]{\fnms{Michael} \snm{Vogt}\corref{}\ead[label=e1]{mv346@cam.ac.uk}}
\runauthor{M. Vogt}
\affiliation{University of Cambridge}
\address[A]{Department of Economics\\
University of Cambridge\\
Sidgwick Avenue\\
Cambridge CB3 9DD\\
United Kingdom\\
\printead{e1}} %adresu isvedimo komanda gale!
\end{aug}

% HISTORY:
\received{\smonth{10} \syear{2011}}
\revised{\smonth{7} \syear{2012}}

% ABSTRACT
%
\begin{abstract}
In this paper, we study nonparametric models allowing for locally
stationary regressors and a regression function that changes smoothly
over time. These models are a natural extension of time series models
with time-varying coefficients. We introduce a kernel-based method to
estimate the time-varying regression function and provide asymptotic
theory for our estimates. Moreover, we show that the main conditions of
the theory are satisfied for a large class of nonlinear autoregressive
processes with a time-varying regression function. Finally, we examine
structured models where the regression function splits up into
time-varying additive components. As will be seen, estimation in these
models does not suffer from the curse of dimensionality.
\end{abstract}

% KEYWORDS
% Pirmas kwd is didziosios raides
%

\end{frontmatter}

%s1 #&#
\section{Introduction}\label{intro}

Classical time series analysis is based on the assumption of
stationarity. However, many time series exhibit a nonstationary
behavior. Examples come from fields as diverse as finance, sound
analysis and neuroscience.

One way to model nonstationary behavior is provided by the theory of
locally stationary processes introduced by Dahlhaus; cf.
\cite{Dahlhaus1996a,Dahlhaus1996b} and~\cite{Dahlhaus1997}. Intuitively
speaking, a process is locally stationary if over short periods of time
(i.e., locally in time) it behaves in an approximately stationary way.
So far, locally stationary models have been mainly considered within a
parametric context. Usually, parametric models are analyzed in which
the coefficients are allowed to change smoothly over time.

There is a considerable amount of papers that deal with time series
models with time-varying coefficients. Dahlhaus et al.
\cite{Dahlhaus1999}, for example, study wavelet estimation in
autoregressive models with time-dependent parameters. Dahlhaus and
Subba Rao~\cite{Dahlhaus2006} analyze a class of ARCH models with
time-varying coefficients. They propose a kernel-based quasi-maximum
likelihood method to estimate the parameter functions; a kernel-based
normalized-least-squares method is suggested by Fryzlewicz et al.
\cite{SubbaRao2008}. Hafner and Linton~\cite{Linton2010} provide
estimation theory for a multivariate GARCH model with a time-varying
unconditional variance. Finally, a~diffusion process with a
time-dependent drift and diffusion function is investigated in Koo and
Linton~\cite{Koo2010}.

In this paper, we introduce a nonparametric framework which can be
regarded as a natural extension of time series models with time-varying
coefficients. In its most general form, the model is given by
%
%e1 #&#
%
\begin{equation}
\label{model} Y_{t,T} = m \biggl(\frac{t}{T},X_{t,T}
\biggr) + \varepsilon_{t,T} \qquad\mbox{for } t = 1,\ldots,T
\end{equation}
with $\ex[\varepsilon_{t,T}|X_{t,T}] = 0$, where $Y_{t,T}$ and
$X_{t,T}$ are random variables of dimension~$1$ and $d$, respectively.
The model variables are assumed to be locally stationary and the
regression function as a whole is allowed to change smoothly over time.
As usual in the literature on locally stationary processes, the
function $m$ does not depend on real time $t$ but rather on rescaled
time $\frac{t}{T}$. This\vspace*{1pt} goes along with the model
variables forming a triangular array instead of a sequence. Throughout
the \hyperref[intro]{Introduction}, we stick to an intuitive concept of
local stationarity. A technically rigorous definition is given in
Section~\ref{section-loc-stat}.

There is a wide range of interesting nonlinear time series models that
fit into the general framework (\ref{model}). An important example is
the nonparametric autoregressive model
%
%e2 #&#
%
\begin{equation}
\label{model-AR} \X_{t,T} = m \biggl(\frac{t}{T},
\X_{t-1,T},\ldots,\X_{t-d,T} \biggr) + \varepsilon_{t,T}
\qquad\mbox{for } t = 1,\ldots,T
\end{equation}
with $\ex[\varepsilon_{t,T} | \X_{t-1,T},\ldots,\X_{t-d,T}]=0$,
which is analyzed in Section~\ref{section-tvNAR}. As will be seen
there, the process defined in (\ref{model-AR}) is locally stationary
and strongly mixing under suitable conditions on the function $m$ and
the error terms $\varepsilon_{t,T}$. Note that independently of the
present work, Kristensen~\cite{Kristensen2011} has developed results
on local stationarity of the process given in (\ref{model-AR}) under a
set of assumptions similar to ours.

In Section~\ref{section-estimation}, we develop estimation theory for
the nonparametric regression function in the general framework (\ref
{model}). As described there, the regression function is estimated by
nonparametric kernel methods. We provide a complete asymptotic theory
for our estimates. In particular, we derive uniform convergence rates
and an asymptotic normality result. To do so, we split up the estimates
into a variance part and a bias part. In order to control the variance
part, we generalize results on uniform convergence rates for kernel
estimates as provided, for example, in Bosq~\cite{Bosq1998}, Masry
\cite{Masry1996} and Hansen~\cite{Hansen2008}. The locally stationary
behavior of the model variables also changes the asymptotic analysis of
the bias part. In particular, it produces an additional bias term which
can be regarded as measuring the deviation from stationarity.

Even though model (\ref{model}) is theoretically interesting, it has
an important drawback. Estimating the time-varying regression function
in (\ref{model}) suffers from an even more severe curse of
dimensionality problem than in the standard, strictly stationary
setting with a time-invariant regression function. The reason is that
in model (\ref{model}), we fit a fully nonparametric function
$m(u,\cdot)$ locally around \textit{each} rescaled time point $u$.
Compared to the standard case, this means that we additionally smooth
in time direction and thus increase the dimensionality of the
estimation problem by one. This makes the procedure even more data
consuming than in the standard setting and thus infeasible in many applications.

In order to countervail this severe curse of dimensionality, we impose
some structural constraints on the regression function in (\ref
{model}). In particular, we consider additive models of the form
%
%e3 #&#
%
\begin{equation}
\label{model-add} Y_{t,T} = \sum_{j=1}^d
m_j \biggl(\frac{t}{T},X_{t,T}^j \biggr) +
\varepsilon_{t,T} \qquad\mbox{for } t = 1,\ldots,T
\end{equation}
with $X_{t,T} = (X_{t,T}^1,\ldots,X_{t,T}^d)$ and $\ex[\varepsilon
_{t,T}|X_{t,T}] = 0$. In Section~\ref{section-add}, we will show that
the component functions of this model can be estimated with
two-dimensional nonparametric convergence rates, no matter how large
the dimension $d$. In order to do so, we extend the smooth backfitting
approach of Mammen et al.~\cite{Mammen1999} to our setting.

%s2 #&#
\section{Local stationarity}\label{section-loc-stat}

Heuristically speaking, a process $\{X_{t,T}\dvtx
t=1,\ldots,\break T\} _{T=1}^{\infty}$ is\vspace*{1pt} locally
stationary if it behaves approximately stationary locally in time. This
intuitive concept can be turned into a rigorous definition in different
ways. One way is to require that locally around each rescaled time
point $u$, the process $\{X_{t,T}\}$ can be approximated by a
stationary process $\{X_t(u)\dvtx t \in\integers\}$ in a stochastic
sense; cf., for example, Dahlhaus and Subba Rao~\cite{Dahlhaus2006}.
This idea also underlies the following definition.
%
%de2.1 #&#
%
\begin{definition}\label{def-loc-stat}
The process $\{X_{t,T}\}$ is locally stationary if for each rescaled
time point $u \in[0,1]$ there exists an associated process $\{X_t(u)\}
$ with the following two properties:
\begin{longlist}
\item$\{X_t(u)\}$ is strictly stationary with density $f_{X_t(u)}$;
\item it holds that
\[
\bigl\| X_{t,T} - X_t(u) \bigr\| \le\biggl( \biggl| \frac{t}{T} - u
\biggr| + \frac{1}{T} \biggr) U_{t,T}(u) \qquad\mbox{a.s.},
\]
where $\{ U_{t,T}(u) \}$ is a process of positive variables satisfying
$\ex[(U_{t,T}(u))^{\rho}] < C$ for some $\rho> 0$ and $C < \infty$
independent of $u$, $t$, and $T$. $\|\cdot\|$ denotes an arbitrary
norm on $\reals^d$.
\end{longlist}
\end{definition}
Since the $\rho$th moments of the variables $U_{t,T}(u)$ are uniformly
bounded, it holds that $U_{t,T}(u) = O_p(1)$. As a consequence of the
above definition, we thus have
\[
\bigl\| X_{t,T} - X_t(u) \bigr\| = O_p \biggl( \biggl|
\frac{t}{T} - u \biggr| + \frac{1}{T} \biggr).
\]
The constant $\rho$ can be regarded as a measure of how well $X_{t,T}$
is approximated by $X_t(u)$: the larger $\rho$ can be chosen, the less
mass is contained in the tails of the distribution of $U_{t,T}(u)$.
Thus, if $\rho$ is large, then the bound $(| \frac{t}{T} - u | +
\frac{1}{T}) U_{t,T}(u)$ will take rather moderate values for most of
the time. In this sense, the bound and thus the approximation of
$X_{t,T}$ by $X_t(u)$ is getting better for larger $\rho$.

%s3 #&#
\section{Locally stationary nonlinear AR models}\label{section-tvNAR}

In this section, we examine a large class of nonlinear autoregressive
processes with a time-varying regression function that fit into the
general framework (\ref{model}). We show that these processes are
locally stationary and strongly mixing under suitable conditions on the
model components. To shorten notation, we repeatedly make use of the
following abbreviation: for any array of variables $\{ Z_{t,T} \}$, we
let $Z_{t,T}^{t-k}:= (Z_{t-k,T},\ldots,Z_{t,T})$ for $k > 0$.

%s3.1 #&#
\subsection{The time-varying nonlinear AR (tvNAR) process}

We call an array $\{ \X_{t,T}\dvtx t \in\integers\}_{T=1}^{\infty}$ a
time-varying nonlinear autoregressive (tvNAR) process if $\X_{t,T}$
evolves according to the equation
%
%e4 #&#
%
\begin{equation}
\label{tvNAR} \X_{t,T} = m \biggl(\frac{t}{T},
\X_{t-1,T}^{t-d} \biggr) + \sigma\biggl(\frac{t}{T},
\X_{t-1,T}^{t-d} \biggr) \varepsilon_t.
\end{equation}
A tvNAR process is thus an autoregressive process of form (\ref
{model-AR}) with errors $\varepsilon_{t,T} = \sigma(\frac{t}{T},\X
_{t-1,T}^{t-d}) \varepsilon_t$. In the above definition, $m(u,y)$ and
$\sigma(u,y)$ are smooth functions of rescaled time $u$ and $y\in
\reals^d$. We stipulate that for $u \le0$, $m(u,y) = m(0,y)$ and
$\sigma(u,y) = \sigma(0,y)$. Analogously, we set $m(u,y) = m(1,y)$
and $\sigma(u,y) = \sigma(1,y)$ for $u \ge1$. Furthermore, the
variables $\varepsilon_t$ are assumed to be i.i.d. with mean zero.
For each $u \in\reals$, we additionally define the associated process
$\{ \X_t(u)\dvtx t \in\integers\}$ by
%
%e5 #&#
%
\begin{equation}
\label{NAR} \X_t(u) = m \bigl(u,\X_{t-1}^{t-d}(u)
\bigr) + \sigma\bigl(u,\X_{t-1}^{t-d}(u) \bigr)
\varepsilon_t,
\end{equation}
where the rescaled time argument of the functions $m$ and $\sigma$ is
fixed at $u$.

As stipulated above, the functions $m$ and $\sigma$ in (\ref{tvNAR}) do
not change over time for $t \le0$. Put differently, $\X_{t,T} = m
(0,\X_{t-1,T}^{t-d} ) + \sigma(0,\X_{t-1,T}^{t-d} ) \varepsilon_t$ for
all $t \le0$. We can thus assume that $\X_{t,T} = \X_t(0)$ for $t
\le0$. Consequently, if there exists a process $\{ \X_t(0)\}$ that
satisfies the system of equations (\ref{NAR}) for $u=0$, then this
immediately implies the existence of a tvNAR process $\{\X_{t,T}\}$
satisfying (\ref{tvNAR}). As will turn out, under appropriate
conditions there exists a strictly stationary solution $\{ \X_t(u)\}$
to (\ref{NAR}) for each $u \in\reals$, in particular for $u = 0$. We
can thus take for granted that the tvNAR process $\{\X_{t,T}\}$ defined
by (\ref{tvNAR}) exists.

Before we turn to the analysis of the tvNAR process, we compare it to
the framework of Zhou and Wu~\cite{Zhou2009} and Zhou~\cite{Zhou2010}.
Their model is given by the equation $Z_{t,T} = G(\frac
{t}{T},\psi_t)$, where $\psi_t =
(\ldots,\varepsilon_{t-1},\varepsilon_t)$ with i.i.d. variables
$\varepsilon_t$ and $G$ is a measurable function. In their theory, the
variables $Z_t(u) = G(u,\psi_t)$ play the role of a stationary
approximation at $u \in[0,1]$. Under suitable assumptions, we can
iterate equation (\ref{NAR}) to obtain that $\X_t(u) = F(u,\psi_t)$ for
some measurable function~$F$. Note, however, that $\X_{t,T} \ne
F(\frac{t}{T},\psi_t)$ in general. This is due to the fact that when
iterating (\ref{NAR}), we use the same functions $m(u,\cdot)$ and
$\sigma(u,\cdot)$ in each step. In contrast to this, different
functions show up in each step when iterating the tvNAR
variables~$\X_{t,T}$. Thus, the relation between the tvNAR process
$\{\X_{t,T}\}$ and the approximations $\{\X_t(u)\}$ is in general
different from that between the processes $\{ Z_{t,T}\}$ and $\{ Z_t(u)
\}$ in the setting of Zhou and Wu.

%s3.2 #&#
\subsection{Assumptions}

We now list some conditions which are sufficient to ensure that the
tvNAR process is locally stationary and strongly mixing. To start with,
the function $m$ is supposed to satisfy the following conditions:
{\renewcommand\thelonglist{(M\arabic{longlist})}
\renewcommand\labellonglist{\thelonglist}
\begin{longlist}
\item\label{M1} $m$ is absolutely bounded by some constant $\um<
\infty$.
\item\label{M2} $m$ is Lipschitz continuous with respect to rescaled
time $u$, that is, there exists a constant $L < \infty$ such that $|
m(u,y) - m(u',y)| \le L |u - u'|$ for all $y\in
\reals^d$.
\item\label{M3} $m$ is continuously differentiable with respect to
$y$. The partial derivatives $\partial_j m(u,y):= \frac{\partial
}{\partial y_j} m(u,y)$ have the property that for some $K_1 < \infty$,
\[
\sup_{u \in\reals, \|y\|_{\infty} > K_1} \bigl|\partial_j m(u,y)\bigr| \le\delta
< 1.
\]
An exact formula for the bound $\delta$ is given in (\ref{delta}) in
Appendix~\ref{appA}.
\end{longlist}}

\noindent The function $\sigma$ is required to fulfill analogous assumptions.
{\renewcommand\thelonglist{($\Sigma$\arabic{longlist})}
\renewcommand\labellonglist{\thelonglist}
\begin{longlist}
\item\label{Sigma1} $\sigma$ is bounded by some constant $\usigma<
\infty$ from above and by some constant $\lsigma> 0$ from below, that
is, $0 < \lsigma\le\sigma(u,y) \le\usigma< \infty$ for all $u$
and~$y$.
\item\label{Sigma2} $\sigma$ is Lipschitz continuous with respect to
rescaled time $u$.
\item\label{Sigma3} $\sigma$ is continuously differentiable with
respect to $y$. The partial derivatives $\partial_j \sigma(u,y):=
\frac{\partial}{\partial y_j} \sigma(u,y)$ have the property that
for some $K_1 < \infty$,\break $|\partial_j \sigma(u,y)| \le\delta< 1$
for all $u \in\reals$ and $\|y\|_{\infty} > K_1$.
\end{longlist}}

\noindent Finally, the error terms are required to have the following properties.
{\renewcommand\thelonglist{(E\arabic{longlist})}
\renewcommand\labellonglist{\thelonglist}
\begin{longlist}
\item\label{E1} The variables $\varepsilon_t$ are i.i.d. with $\ex
[\varepsilon_t] = 0$ and $\ex|\varepsilon_t|^{1+\eta} < \infty$
for some \mbox{$\eta> 0$}. Moreover, they have an everywhere positive and
continuous density~$f_{\varepsilon}$.
\item\label{E2} The density $f_{\varepsilon}$ is bounded and
Lipschitz, that is, there exists a constant $L < \infty$ such that $|
f_{\varepsilon}(z) - f_{\varepsilon}(z')| \le L |z -
z'|$ for all $z, z' \in\reals$.
\end{longlist}}

\noindent To show that the tvNAR process is strongly mixing, we additionally need
the following condition on the density of the error terms:
{\renewcommand\thelonglist{(E\arabic{longlist})}
\renewcommand\labellonglist{\thelonglist}
\begin{longlist}
\setcounter{longlist}{2}
\item\label{E3} Let $d_0$, $d_1$ be any constants with $0 \le d_0 \le
D_0 < \infty$ and $|d_1| \le D_1 < \infty$. The density
$f_{\varepsilon}$ fulfills the condition
\[
\int_{\reals} \bigl| f_{\varepsilon}\bigl([1 + d_0]z +
d_1\bigr) - f_{\varepsilon
}(z)\bigr| \,dz \le C_{D_0,D_1} \bigl(
d_0 + |d_1| \bigr)
\]
with $C_{D_0,D_1} < \infty$ only depending on the bounds $D_0$ and $D_1$.
\end{longlist}}

We shortly give some remarks on the above conditions:

\begin{longlist}[(iii)]
\item[(i)] Our set of assumptions can be regarded as a strengthening of the
assumptions needed to show geometric ergodicity of nonlinear AR
processes of the form $\X_t = m(\X_{t-1}^{t-d}) + \sigma(\X
_{t-1}^{t-d}) \varepsilon_t$. The main assumption in this context
requires the functions $m$ and $\sigma$ not to grow too fast outside a
large bounded set. More precisely, it requires them to be dominated by
linear functions with sufficiently small slopes; cf. Tj{\o}stheim
\cite{Tjostheim1990}, Bhattacharya and Lee~\cite{Bhattacharya1995},
An and Huang~\cite{An1996} or Chen and Chen~\cite{Chen2000}, among
others.~\ref{M3} and~\ref{Sigma3} are very close in
spirit to this kind of assumption. They restrict the growth of $m$ and
$\sigma$ by requiring the derivatives of these functions to be small
outside a large bounded set.

\item[(ii)] If we replace~\ref{M3} and~\ref{Sigma3} with the
stronger assumption that the partial derivatives $|\partial_j m(u,y)|$
and $|\partial_j \sigma(u,y)|$ are globally bounded by some
sufficiently small number $\delta< 1$, then some straightforward
modifications allow us to dispense with the boundedness assumptions
\ref{M1} and~\ref{Sigma1} in the local stationarity and
mixing proofs.

\item[(iii)] Condition~\ref{M3} implies that the derivatives $\partial_j
m(u,y)$ are absolutely bounded. Hence, there exists a constant $\Delta
< \infty$ such that $|\partial_j m(u,y)| \le\Delta$ for all $u \in
\reals$ and $y\in\reals^d$. Similarly,~\ref{Sigma3}
implies that the derivatives $\partial_j \sigma(u,y)$ are absolutely
bounded by some constant $\Delta< \infty$.

\item[(iv)] As already noted,~\ref{E3} is only needed to prove that the
tvNAR process is strongly mixing. It is, for example, fulfilled for the
class of bounded densities $f_{\varepsilon}$ whose first derivative
$f_{\varepsilon}'$ is bounded, satisfies $\int|z
f_{\varepsilon}'(z)| \,dz < \infty$ and declines monotonically
to zero for values $|z| > C$ for some constant $C > 0$; see also
Section~3 in Fryzlewicz and Subba Rao~\cite{SubbaRao2010} who work
with assumptions closely related to~\ref{E3}.
\end{longlist}

%s3.3 #&#
\subsection{Properties of the tvNAR process}

We now show that the tvNAR process is locally stationary and strongly
mixing under the assumptions listed above. In addition, we will see
that the auxiliary processes $\{\X_t(u)\}$ have densities that vary
smoothly over rescaled time $u$. As will turn out, these three
properties are central for the estimation theory developed in Sections
\ref{section-estimation} and~\ref{section-add}.

The first theorem summarizes some properties of the tvNAR process and
of the auxiliary processes $\{\X_t(u)\}$ that are needed to prove the
main results.

%th3.1 #&#
%
\begin{theorem}\label{theo-stat}
Let~\ref{M1}--\ref{M3},~\ref{Sigma1}--\ref{Sigma3} and
\ref{E1} be fulfilled. Then:

{\renewcommand\thelonglist{(\roman{longlist})}
\renewcommand\labellonglist{\thelonglist}
\begin{longlist}
\item\label{theo-stat-1} for each $u \in\reals$, the process $\{\X
_t(u), t \in\integers\}$ has a strictly stationary solution with
$\varepsilon_t$ independent of $\X_{t-k}(u)$ for $k > 0$;
\item\label{theo-stat-2} the variables $\X_{t-1}^{t-d}(u)$ have a
density $f_{\X_{t-1}^{t-d}(u)}$ w.r.t. Lebesgue measure;
\item\label{theo-stat-3} the variables $\X_{t-1,T}^{t-d}$ have
densities $f_{\X_{t-1,T}^{t-d}}$ w.r.t. Lebesgue measure.
\end{longlist}}
\end{theorem}

The next result states that $\{\X_{t,T}\}$ can be locally approximated
by $\{\X_t(u)\}$. Together with Theorem~\ref{theo-stat}, it shows
that the tvNAR process $\{\X_{t,T}\}$ is locally stationary in the
sense of Definition~\ref{def-loc-stat}.
%
%th3.2 #&#
%
\begin{theorem}\label{theo-loc-stat}
Let~\ref{M1}--\ref{M3},~\ref{Sigma1}--\ref{Sigma3} and
\ref{E1} be fulfilled. Then
%
%e6 #&#
%
\begin{equation}
\label{loc-stat-inequ} \bigl| \X_{t,T} - \X_t(u) \bigr| \le\biggl( \biggl|
\frac{t}{T} - u \biggr| + \frac{1}{T} \biggr) U_{t,T}(u)
\qquad\mbox{a.s.},
\end{equation}
where the variables $U_{t,T}(u)$ have the property that $\ex
[(U_{t,T}(u))^{\rho}] < C$ for some $\rho> 0$ and $C < \infty$
independent of $u$, $t$ and $T$.
\end{theorem}
To get an idea of the proof of Theorem~\ref{theo-loc-stat}, consider
the model $\X_{t,T} = m(\frac{t}{T},\X_{t-1,T}) + \varepsilon_t$
for a moment. Our arguments are based on a backward expansion of the
difference $\X_{t,T} - \X_t(u)$. Exploiting the smoothness conditions
of~\ref{M2} and~\ref{M3} together with the boundedness of $m$,
we obtain that
\[
\bigl| \X_{t,T} - \X_t(u) \bigr| \le C \sum
_{r=0}^{n-1} \prod_{k=1}^r
\bigl| \partial m(u,\xi_{t-k}) \bigr| \biggl( \biggl| \frac{t}{T} - u \biggr| +
\frac{r}{T} \biggr) + C \prod_{k=1}^n
\bigl| \partial m(u,\xi_{t-k}) \bigr|,
\]
where $\partial m(u,y)$ is the derivative of $m(u,y)$ with respect to
$y$ and $\xi_{t-k}$ is an intermediate point between $\X_{t-k,T}$ and
$\X_{t-k}(u)$. To prove (\ref{loc-stat-inequ}), we have to show that
the product $\prod_{k=1}^n | \partial m(u,\xi_{t-k}) |$ is
contracting in some stochastic sense as $n$ tends to infinity. The
heuristic idea behind the proof is the following: using conditions
\ref{M1} and~\ref{E1}, we can show that at least a certain
fraction of the terms $\xi_{t-1},\ldots,\xi_{t-n}$ take a value in
the region $\{ y\dvtx|y| > K_1 \}$ as $n$ grows large.\vadjust{\goodbreak} Since the
derivative $|\partial m|$ is small in this region according to
\ref{M3}, this ensures that at least a certain fraction of the elements in
the product $\prod_{k=1}^n | \partial m(u,\xi_{t-k}) |$ are
small in value. This prevents the product from exploding and makes it
contract to zero as $n$ goes to infinity.

Next, we come to a result which shows that the densities of the
approximating variables $\X_{t-1}^{t-d}(u)$ change smoothly over time.
%
%th3.3 #&#
%
\begin{theorem}\label{theo-smooth-dens}
Let $f(u,y):= f_{\X_{t-1}^{t-d}(u)}(y)$ be the density of $\X
_{t-1}^{t-d}(u)$ at \mbox{$y\in\reals^d$}. If~\ref{M1}--\ref{M3},
\ref{Sigma1}--\ref{Sigma3} and~\ref{E1},~\ref{E2} are fulfilled, then
\[
\bigl|f(u,y) - f(v,y)\bigr| \le C_{y} |u - v|^p
\]
with some constant $0 < p < 1$ and $C_{y} < \infty$ continuously
depending on $y$.
\end{theorem}

We finally characterize the mixing behavior of the tvNAR process. To do
so, we first give a quick reminder of the definitions of an $\alpha$-
and $\beta$-mixing array.
Let $(\Omega,\mathcal{A},\pr)$ be a probability space, and let
$\mathcal{B}$ and $\mathcal{C}$ be subfields of~$\mathcal{A}$. Define
\begin{eqnarray*}
\alpha(\mathcal{B},\mathcal{C}) & = & \sup_{B \in\mathcal{B}, C \in
\mathcal{C}} \bigl| \pr(B \cap C) -
\pr(B) \pr(C) \bigr|,
\\
\beta(\mathcal{B},\mathcal{C}) & = & \ex\sup_{C \in\mathcal{C}} \bigl| \pr(C)
- \pr(C|
\mathcal{B}) \bigr|.
\end{eqnarray*}
Moreover, for an array $\{ Z_{t,T}\dvtx1 \le t \le T \}$, define the
coefficients
%
%e7 #&#
%e8 #&#
%
\begin{eqnarray}
\label{alpha-coeff}
\alpha(k) & = & \sup_{t,T\dvtx1 \le t \le T-k} \alpha\bigl( \sigma
(Z_{s,T},1 \le
s \le t),\sigma(Z_{s,T},t+k \le s \le T) \bigr),
\\
\label{beta-coeff}
\beta(k) & = & \sup_{t,T\dvtx1 \le t \le T-k} \beta\bigl( \sigma(Z_{s,T},1
\le s
\le t),\sigma(Z_{s,T},t+k \le s \le T) \bigr),
\end{eqnarray}
where $\sigma(Z)$ is the $\sigma$-field generated by $Z$. The array
$\{Z_{t,T}\}$ is said to be $\alpha$-mixing (or strongly mixing) if
$\alpha(k) \rightarrow0$ as $k \rightarrow\infty$. Similarly, it is
called $\beta$-mixing if $\beta(k) \rightarrow0$. Note that $\beta
$-mixing implies $\alpha$-mixing. The final result of this section
shows that the tvNAR process is $\beta$-mixing with coefficients that
converge exponentially fast to zero.
%
%th3.4 #&#
%
\begin{theorem}\label{theo-mixing}
If~\ref{M1}--\ref{M3},~\ref{Sigma1}--\ref{Sigma3} and
\ref{E1}--\ref{E3} are fulfilled, then the tvNAR process $\{\X_{t,T}\}$
is geometrically $\beta$-mixing, that is, there exist positive
constants $\gamma< 1$ and $C < \infty$ such that $\beta(k) \le C
\gamma^k$.
\end{theorem}
The strategy of the proof is as follows: the (conditional)
probabilities that show up in the definition of the $\beta
$-coefficient in (\ref{beta-coeff}) can be written in terms of the
functions $m$, $\sigma$ and the error density $f_{\varepsilon}$. To
do so, we derive recursive expressions of the model variables $\X
_{t,T}$ and of certain conditional densities of $\X_{t,T}$. Rewriting
the $\beta$-coefficient with the help of these expressions allows us
to derive an appropriate\vadjust{\goodbreak} bound for it. The overall strategy is thus
similar to that of Fryzlewicz and Subba Rao~\cite{SubbaRao2010} who
also derive bounds of mixing coefficients in terms of conditional
densities. The specific steps of the proof, however, are quite
different. The details together with the proofs of the other theorems
can be found in Appendix~\ref{appA}.\vspace*{-2pt}

%s4 #&#
\section{Kernel estimation}\label{section-estimation}

In this section, we consider kernel estimation in the general model
(\ref{model}),
\[
Y_{t,T} = m \biggl(\frac{t}{T},X_{t,T} \biggr) +
\varepsilon_{t,T} \qquad\mbox{for } t=1,\ldots,T
\]
with $\ex[\varepsilon_{t,T}|X_{t,T}] = 0$. Note that $m(\frac
{t}{T},\cdot)$ is the conditional mean function in model (\ref
{model}) at the time point $t$. The function $m$ is thus identified
almost surely on the grid of points $\frac{t}{T}$ for $t=1,\ldots,T$.
These points form a dense subset of the unit interval as the sample
size grows to infinity. As a consequence, $m$ is identified almost
surely at all rescaled time points $u \in[0,1]$ if it is continuous in
time direction (which we will assume in what follows).\vspace*{-2pt}

%s4.1 #&#
\subsection{Estimation procedure}

We restrict attention to Nadaraya--Watson (NW) estimation. It is
straightforward to extend the theory to local linear (or more generally
local polynomial) estimation. The NW estimator of model (\ref{model})
is given by
%
%e9 #&#
%
\begin{equation}
\label{nw} \hat{m}(u,x) = \frac{\sum_{t=1}^T K_h (u -
{t}/{T} ) \prod_{j=1}^d K_h (x^j-X_{t,T}^j )
Y_{t,T}}{\sum_{t=1}^T K_h (u-{t}/{T} ) \prod_{j=1}^d K_h
(x^j-X_{t,T}^j )}.
\end{equation}
Here and in what follows, we write $X_{t,T} =
(X_{t,T}^1,\ldots,X_{t,T}^d)$ and $x = (x^1,\ldots,x^d)$ for any vector
$x \in\reals^d$, that is, we use subscripts to indicate the time point
of observation and superscripts to denote the components of the vector.
$K$ denotes a one-dimensional kernel function and we use the notation
$K_h(v) = K(\frac{v}{h})$. For convenience, we work with a product
kernel and assume that the bandwidth $h$ is the same in each direction.
Our results can, however, be easily modified to allow for nonproduct
kernels and different bandwidths.

The estimate defined in (\ref{nw}) differs from the NW estimator in
the standard strictly stationary setting in that there is an additional
kernel in time direction. We thus do not only smooth in the direction
of the covariates $X_{t,T}$ but also in the time direction. This takes
into account that the regression function is varying over time. In what
follows, we derive the asymptotic properties of our NW estimate. The
proofs are given in Appendix~\ref{appB}.\vspace*{-2pt}

%s4.2 #&#
\subsection{Assumptions}

The following three conditions are central to our results:
{\renewcommand\thelonglist{(C\arabic{longlist})}
\renewcommand\labellonglist{\thelonglist}
\begin{longlist}
\item\label{C1} The process $\{X_{t,T}\}$ is locally stationary in
the sense of Definition~\ref{def-loc-stat}.
Thus, for each time point $u \in[0,1]$, there exists a strictly
stationary process $\{X_t(u)\}$ having the property that $\|X_{t,T} -
X_t(u)\| \le(|\frac{t}{T}-u| + \frac{1}{T}) U_{t,T}(u)$ a.s. with
$\ex[(U_{t,T}(u))^{\rho}] \le C$ for some $\rho> 0$.\vadjust{\goodbreak}
\item\label{C2} The densities $f(u,x):= f_{X_t(u)}(x)$ of the
variables $X_t(u)$ are smooth in $u$. In particular, $f(u,x)$ is
differentiable w.r.t. $u$ for each $x \in\reals^d$, and the
derivative $\partial_0 f(u,x):= \frac{\partial}{\partial u} f(u,x)$
is continuous.
\item\label{C3} The array $\{ X_{t,T}, \varepsilon_{t,T} \}$ is
$\alpha$-mixing.
\end{longlist}}

\noindent As seen in Section~\ref{section-tvNAR}, these three
conditions are essentially fulfilled for the tvNAR process:~\ref{C1}
and~\ref{C3} follow immediately from Theorems~\ref{theo-loc-stat} and~\ref {theo-mixing}.
Moreover, Theorem~\ref{theo-smooth-dens} shows that
the tvNAR process satisfies a weakened version of~\ref{C2} which
requires the densities $f_{X_t(u)}$ to be continuous rather than
differentiable in time direction. Note that we could do with this
weakened version of~\ref{C2}, however at the cost of getting slower
convergence rates for the bias part of the NW estimate.

In addition to the above three assumptions, we impose the following
regularity conditions:
{\renewcommand\thelonglist{(C\arabic{longlist})}
\renewcommand\labellonglist{\thelonglist}
\begin{longlist}
\setcounter{longlist}{3}
\item\label{C4} $f(u,x)$ is partially differentiable w.r.t. $x$ for
each $u \in[0,1]$. The derivatives $\partial_j f(u,x):= \frac
{\partial}{\partial x^j} f(u,x)$ are continuous for $j=1,\ldots,d$.
\item\label{C5} $m(u,x)$ is twice continuously partially
differentiable with first derivatives $\partial_j m(u,x)$ and second
derivatives $\partial_{ij}^2 m(u,x)$ for $i,j=0,\ldots,d$.
\item\label{C6} The kernel $K$ is symmetric about zero, bounded and
has compact support, that is, $K(v) = 0$ for all $|v| > C_1$ with some
$C_1 < \infty$. Furthermore, $K$ is Lipschitz, that is, $| K(v) -
K(v') | \le L | v - v' |$ for some $L < \infty$ and
all $v,v' \in\reals$.
\end{longlist}}

\noindent Finally, note that throughout the paper the bandwidth $h$ is assumed to
converge to zero at least at polynomial rate, that is, there exists a
small $\xi> 0$ such that $h \le C T^{-\xi}$ for some constant $C >
0$.\vspace*{-2pt}

%s4.3 #&#
\subsection{Uniform convergence rates for kernel averages}
\label{subsection-conv-rates-1}

As a first step in the analysis of the NW estimate (\ref{nw}), we
examine kernel averages of the general form
%
%e10 #&#
%
\begin{equation}
\label{kernel-average-1} \hat{\psi}(u,x) = \frac{1}{T h^{d+1}} \sum
_{t=1}^T K_h \biggl(u -
\frac{t}{T} \biggr) \prod_{j=1}^d
K_h\bigl(x^j - X_{t,T}^j\bigr)
W_{t,T}
\end{equation}
with $\{ W_{t,T} \}$ being an array of one-dimensional random
variables. A wide range of kernel-based estimators, including the NW
estimator defined in (\ref{nw}), can be written as functions of
averages of the above form. The asymptotic behavior of such averages is
thus of wider interest. For this reason, we investigate the properties
of these averages for a general array of variables $\{ W_{t,T} \}$.
Later on we will employ the results with $W_{t,T} = 1$ and $W_{t,T} =
\varepsilon_{t,T}$.

We now derive the uniform convergence rate of $\hat{\psi}(u,x) - \ex
\hat{\psi}(u,x)$. To do so, we make the following assumptions on the
components in (\ref{kernel-average-1}):
{\renewcommand\thelonglist{(K\arabic{longlist})}
\renewcommand\labellonglist{\thelonglist}
\begin{longlist}
\item\label{K0} It holds that $\ex|W_{t,T}|^s \le C$ for some $s >
2$ and $C < \infty$.
\item\label{K1} The array $\{X_{t,T},W_{t,T}\}$ is $\alpha$-mixing.
The mixing coefficients $\alpha$ have the property that $\alpha(k)
\le A k^{-\beta}$ for some $A < \infty$ and $\beta> \frac{2s-2}{s-2}$.\vadjust{\goodbreak}
\item\label{K2} Let $f_{X_{t,T}}$ and $f_{X_{t,T},X_{t+l,T}}$ be the
densities of $X_{t,T}$ and\break $(X_{t,T},X_{t+l,T})$, respectively. For any
compact set $S \subseteq\reals^d$, there exists a \mbox{constant} $C = C(S)$
such that $\sup_{t,T} \sup_{x \in S} f_{X_{t,T}}(x) \le C$ and\break
\mbox{$\sup_{t,T} \sup_{x \in S} \ex[ |W_{t,T}|^s | X_{t,T}=x
] f_{X_{t,T}}(x) \le C$}.
Moreover, there exists\break a \mbox{natural} number $l^* < \infty$ such that for
all $l \ge l^*$,\break
$\sup_{t,T} \sup_{x,x' \in S} \ex[ |W_{t,T}|
|W_{t+l,T}| | X_{t,T}=x, X_{t+l,T} = x' ]
f_{X_{t,T},X_{t+l,T}}(x,x') \le C$.
\end{longlist}}

The next theorem generalizes uniform convergence results of Hansen
\cite{Hansen2008} for the strictly stationary case to our setting. See
Kristensen~\cite{Kristensen2009} for related results.
%
%th4.1 #&#
%
\begin{theorem}\label{theo-stochastic-part}
Assume that~\ref{K0}--\ref{K2} are satisfied with
%
%e11 #&#
%
\begin{equation}
\label{beta}
\beta> \frac{2 + s(1 + (d+1))}{s-2}
\end{equation}
and that the kernel $K$ fulfills~\ref{C6}. In addition, let the
bandwidth satisfy
%
%e12 #&#
%
\begin{equation}
\label{h} \frac{\phi_T \log T}{T^{\theta} h^{d+1}} = o(1)
\end{equation}
with $\phi_T$ slowly diverging to infinity (e.g., $\phi_T = \log
\log T$) and
%
%e13 #&#
%
\begin{equation}
\label{theta} \theta= \frac{\beta(1-{2}/{s}) - {2}/{s} - 1 -
(d+1)}{\beta+ 3 - (d+1)}.
\end{equation}
Finally, let $S$ be a compact subset of $\reals^d$. Then it holds that
%
%e14 #&#
%
\begin{equation}
\sup_{u \in[0,1], x \in S} \bigl| \hat{\psi}(u,x) - \ex\hat{\psi}(u,x) \bigr| = O_p
\biggl(\sqrt{\frac{\log T}{T h^{d+1}}} \biggr).
\end{equation}
\end{theorem}
The convergence rate in the above theorem is identical to the rate
obtained for a $(d+1)$-dimensional nonparametric estimation problem in
the standard strictly stationary setting. This reflects the fact that
additionally smoothing in time direction, we essentially have a
$(d+1)$-dimensional problem in our case. Moreover, note that with (\ref
{beta}) and (\ref{theta}), we can compute that $\theta\in(0,1-\frac
{2}{s}]$. In particular, $\theta= 1-\frac{2}{s}$ if the mixing
coefficients decay exponentially fast to zero, that is, if $\beta=
\infty$. Restriction (\ref{h}) on the bandwidth is thus a
strengthening of the usual condition that $T h^{d+1} \rightarrow\infty$.

%s4.4 #&#
\subsection{Uniform convergence rates for NW estimates}
\label{subsection-conv-rates-3}

The next theorem characterizes the uniform convergence behavior of our
NW estimate.
%
%th4.2 #&#
%
\begin{theorem}\label{theo-nw} Assume that~\ref{C1}--\ref{C6}
hold and that~\ref{K0}--\ref{K2} are fulfilled both for
$W_{t,T} = 1$ and $W_{t,T} = \varepsilon_{t,T}$. Let $\beta$ satisfy
(\ref{beta}) and suppose that $\inf_{u \in[0,1], x \in S} f(u,x) >
0$. Moreover, assume that the bandwidth $h$ satisfies
%
%e15 #&#
%
\begin{equation}
\frac{\phi_T \log T}{T^{\theta} h^{d+1}} = o(1) \quad\mbox{and}\quad
\frac{1}{T^r
h^{d+r}} = o(1)
\end{equation}
with $\theta$ given in (\ref{theta}), $\phi_T = \log\log T$, $r =
\min\{\rho,1\}$ and $\rho$ introduced in~\ref{C1}. Defining $I_h
= [C_1h, 1 - C_1h]$, it then holds that
%
%e16 #&#
%
\begin{equation}
\label{rate-nw} \sup_{u \in I_h, x \in S} \bigl| \hat{m}(u,x) - m(u,x)
\bigr| = O_p \biggl( \sqrt{\frac{\log T}{T h^{d+1}}} + \frac{1}{T^r
h^d} + h^2 \biggr).
\end{equation}
\end{theorem}
To derive the above result, we decompose the difference $\hat{m}(u,x)
- m(u,x)$ into a stochastic part and a bias part. Using Theorem \ref
{theo-stochastic-part}, the stochastic part can be shown to be of the
order $O_p(\sqrt{\log T / T h^{d+1}})$. The bias term splits up into
two parts, a standard component of the order $O(h^2)$ and a nonstandard
component of the order $O(T^{-r} h^{-d})$. The latter component results
from replacing the variables $X_{t,T}$ by $X_t(\frac{t}{T})$ in the
bias term. It thus captures how far these variables are from their
stationary approximations $X_t(\frac{t}{T})$. Put differently, it
measures the deviation from stationarity. As will be seen in Appendix
\ref{appB}, handling this nonstationarity bias requires techniques substantially
different from those needed to treat the bias term in a strictly
stationary setting.

Note that the additional nonstationarity bias converges faster to zero
for larger $r = \min\{\rho,1\}$. This makes perfect sense if we
recall from Section~\ref{section-loc-stat} that $r$ measures how well
$X_{t,T}$ is locally approximated by $X_t(\frac{t}{T})$: the larger
$r$, the smaller the deviation of $X_{t,T}$ from its stationary
approximation and thus the smaller the additional nonstationarity bias.

%s4.5 #&#
\subsection{Asymptotic normality}\label{subsection-normality}

We conclude the asymptotic analysis of our NW estimate with a result on
asymptotic normality.
%
%th4.3 #&#
%
\begin{theorem}\label{theo-normality}
Assume that~\ref{C1}--\ref{C6} hold and that~\ref{K0}--\ref{K2} are
fulfilled both for $W_{t,T} = 1$ and $W_{t,T} = \varepsilon_{t,T}$. Let
$\beta\ge4$ and $T^r h^{d+2} \rightarrow \infty$ with $r =
\min\{\rho,1\}$. Moreover, suppose that $f(u,x) > 0$ and that
$\sigma^2(\frac{t}{T},x):= \ex[\varepsilon_{t,T}^2 | X_{t,T} = x]$ is
continuous. Finally, let $r > \frac{d+2}{d+5}$ to ensure that the
bandwidth $h$ can be chosen to satisfy $T h^{d+5} \rightarrow c_h$ for
a constant $c_h$. Then
%
%e17 #&#
%
\begin{equation}
\sqrt{T h^{d+1}} \bigl( \hat{m}(u,x) - m(u,x) \bigr) \convd\normal
(B_{u,x},V_{u,x}),
\end{equation}
where $B_{u,x} = \sqrt{c_h} \frac{\kappa_2}{2} \sum_{i=0}^d [ 2 \,\partial
_i m(u,x) \,\partial_i f(u,x) + \partial_{i,i}^2
m(u,x) f(u,x) ] / f(u,x)$ and $V_{u,x} = \kappa_0^{d+1} \sigma^2(u,x)
/ f(u,x)$ with $\kappa_0 = \int K^2(\varphi) \,d\varphi$ and $\kappa_2 =
\int\varphi^2 K(\varphi) \,d\varphi$.
\end{theorem}
The above theorem parallels the asymptotic normality result for the
standard strictly stationary setting. In particular, the bias and
variance expressions $B_{u,x}$ and $V_{u,x}$ are very\vspace*{1pt}
similar to those from the standard case. By requiring that $T^r h^{d+2}
\rightarrow \infty$, we make sure that the additional nonstationarity
bias is asymptotically negligible.

%s5 #&#
\section{Locally stationary additive models}\label{section-add}

We now put some structural constraints on the regression function $m$
in model (\ref{model}). In particular, we assume that for all rescaled
time points $u \in[0,1]$ and all points $x$ in a compact subset of
$\reals^d$, say $[0,1]^d$, the regression function can be split up into
additive components according to $m(u,x) = m_0(u) + \sum_{j=1}^d
m_j(u,x^j)$. This means that for $x \in[0,1]^d$, we have the additive
regression model
%
%e18 #&#
%
\begin{equation}
\label{additive-model} \ex[Y_{t,T}|X_{t,T} = x] =
m_0 \biggl(\frac{t}{T} \biggr) + \sum
_{j=1}^d m_j \biggl(
\frac{t}{T},x^j \biggr).
\end{equation}

To identify the component functions of model (\ref{additive-model})
within the unit cube $[0,1]^d$, we impose\vspace*{1pt} the condition that $\int
m_j(u,x^j) p_j(u,x^j) \,dx^j = 0$ for all $j = 1,\ldots,d$ and all
rescaled time points $u \in[0,1]$. Here, the functions $p_j(u,x^j) =
\int p(u,x) \,dx^{-j}$ are the marginals of the density
\[
p(u,x) = \frac{I(x \in[0,1]^d) f(u,x)}{\pr(X_0(u) \in[0,1]^d)},
\]
where as before $f(u,\cdot)$ is the density of the strictly stationary
process $\{X_t(u)\}$. Note that this normalization of the component
functions varies over time in the sense that for each rescaled time
point $u$, we integrate with respect to a different density.

To estimate the functions $m_0,\ldots,m_d$, we adapt the smooth
backfitting technique of Mammen et al.~\cite{Mammen1999} to our
setting. To do so, we first introduce the auxiliary estimates
\begin{eqnarray*}
\hat{p}(u,x) & = & \frac{1}{\Obs} \sum_{t=1}^T
I\bigl(X_{t,T} \in[0,1]^d\bigr) K_h \biggl(u,
\frac{t}{T} \biggr) \prod_{j=1}^d
K_h \bigl(x^j,X_{t,T}^j \bigr),
\\
\hat{m}(u,x) & = & \frac{1}{\Obs} \sum_{t=1}^T
I\bigl(X_{t,T} \in[0,1]^d\bigr) K_h \biggl(u,
\frac{t}{T} \biggr) \prod_{j=1}^d
K_h \bigl(x^j,X_{t,T}^j \bigr)
Y_{t,T} / \hat{p}(u,x).
\end{eqnarray*}
$\hat{p}(u,x)$ is a kernel estimate of the density $p(u,x)$, and $\hat
{m}(u,x)$ is a $(d+1)$-dimensional NW smoother that estimates $m(u,x)$
for $x \in[0,1]^d$. In the above definitions,
\[
\Obs= \sum_{t=1}^T K_h
\biggl(u,\frac{t}{T} \biggr) I\bigl(X_{t,T} \in[0,1]^d
\bigr)
\]
is the number of observations in the unit cube $[0,1]^d$, where only
time points close to $u$ are taken into account, and
\[
K_h(v,w) = I\bigl(v,w \in[0,1]\bigr) \frac{K_h(v-w)}{\int_0^1 K_h(s-w)
\,ds}
\]
is a modified kernel weight. This weight has the property that
\mbox{$\int_0^1\! K_h(v, w) \,dv\! =\! 1$} for all $w \in[0,1]$, which is needed to
derive the asymptotic properties of the backfitting estimates.

Given the smoothers $\hat{p}$ and $\hat{m}$, we define the smooth
backfitting estimates $\tilde{m}_0(u)$,
$\tilde{m}_1(u,\cdot),\ldots,\tilde{m}_d(u,\cdot)$ of the functions
$m_0(u)$, $m_1(u,\cdot),\ldots,m_d(u,\cdot)$ at the time point $u
\in[0,1]$ as the minimizers of the criterion
%
%e19 #&#
%
\begin{equation}
\label{proj} \int\Biggl( \hat{m}(u,w) - g_0 - \sum
_{j=1}^d g_j\bigl(w^j\bigr)
\Biggr)^2 \hat{p}(u,w) \,dw,
\end{equation}
where the minimization runs over all additive functions $g(x) = g_0 +
g_1(x^1) + \cdots+ g_d(x^d)$ whose components are normalized to satisfy
$\int g_j(w^j) \hat{p}_j(u,\break w^j) \,dw^j = 0$ for $j = 1,\ldots,d$. Here,
$\hat{p}_j(u,x^j) = \int\hat{p}(u,x) \,dx^{-j}$ is the marginal of the
kernel density $\hat{p}(u,\cdot)$ at the point $x^j$.

According to (\ref{proj}), the backfitting estimate $\tilde
{m}(u,\cdot) = \tilde{m}_0(u) + \sum_{j=1}^d \tilde
{m}_j(u,\cdot)$ is an $L_2$-projection of the full-dimensional NW
estimate $\hat{m}(u,\cdot)$ onto the subspace of additive functions,
where the projection is done with respect to the density estimate $\hat
{p}(u,\cdot)$. Note that (\ref{proj}) is a $d$-dimensional projection
problem. In particular, rescaled time does not enter as an additional
dimension. The projection is rather done separately for each time point
$u \in[0,1]$. We thus fit a smooth backfitting estimate to the data
separately around each point in time $u$.

By differentiation, we can show that the minimizer of (\ref{proj}) is
characterized by the system of integral equations
%
%e20 #&#
%
\begin{equation}
\label{int-equ}\qquad \tilde{m}_j\bigl(u,x^j\bigr) =
\hat{m}_j\bigl(u,x^j\bigr) - \sum
_{k \ne j} \int\tilde{m}_k\bigl(u,x^k
\bigr) \frac{\hat{p}_{j,k}(u,x^j,x^k)}{\hat
{p}_j(u,x^j)} \,dx^k - \tilde{m}_0(u)
\end{equation}
together with $\int\tilde{m}_j(u,w^j) \hat{p}_j(u,w^j) \,dw^j = 0$ for
$j=1,\ldots,d$. Here, $\hat{p}_j$ and $\hat{p}_{j,k}$ are kernel
density estimates, and $\hat{m}_j$ is a NW smoother defined as
\begin{eqnarray*}
\hat{p}_j\bigl(u,x^j\bigr) & = & \frac{1}{\Obs} \sum
_{t=1}^T I\bigl(X_{t,T}
\in[0,1]^d\bigr) K_h \biggl(u,\frac{t}{T} \biggr)
K_h\bigl(x^j,X_{t,T}^j\bigr),
\\
\hat{p}_{j,k}\bigl(u,x^j,x^k\bigr) & = &
\frac{1}{\Obs} \sum_{t=1}^T I
\bigl(X_{t,T} \in[0,1]^d\bigr) K_h \biggl(u,
\frac{t}{T} \biggr)\\
&&\qquad\hspace*{22pt}{}\times K_h\bigl(x^j,X_{t,T}^j
\bigr) K_h\bigl(x^k,X_{t,T}^k\bigr),
\\
\hat{m}_j\bigl(u,x^j\bigr) & = & \frac{1}{\Obs} \sum
_{t=1}^T I\bigl(X_{t,T}
\in[0,1]^d\bigr) K_h \biggl(u,\frac{t}{T} \biggr)\\
&&\qquad\hspace*{22pt}{}\times
K_h\bigl(x^j,X_{t,T}^j\bigr)
Y_{t,T} / \hat{p}_j\bigl(u,x^j\bigr).
\end{eqnarray*}
Moreover, the estimate $\tilde{m}_0(u)$ of the model constant at time
point $u$ is given by $\tilde{m}_0(u) = \Obs^{-1} \sum_{t=1}^T
I(X_{t,T} \in[0,1]^d) K_h(u,\frac{t}{T}) Y_{t,T}$.

We next summarize the assumptions needed to derive the asymptotic
properties of the smooth backfitting estimates. First of all, the
conditions of Section~\ref{section-estimation} must be satisfied for
the kernel estimates that show up in the system of integral
equations~(\ref{int-equ}). This is ensured by the following
assumption.\looseness=-1
{\renewcommand\thelonglist{(Add\arabic{longlist})}
\renewcommand\labellonglist{\thelonglist}
\begin{longlist}
\item\label{add1} Conditions~\ref{C1}--\ref{C6} are fulfilled together
with~\ref{K0}--\ref{K2} for $W_{t,T}=1$ and
$W_{t,T}=\varepsilon_{t,T}$. The parameter $\beta$ satisfies the
inequality $\beta> \max\{ 4,\frac{2 + 3s}{s-2} \}$ and $\inf_{u \in
[0,1], x \in[0,1]^d} f(u,x) > 0$.
\end{longlist}}

\noindent In addition to~\ref{add1}, we need some restrictions on the
admissible bandwidth. For convenience, we stipulate somewhat stronger
conditions than in Section~\ref{section-estimation} to get rid of the
additional nonstationarity bias from the very beginning.
{\renewcommand\thelonglist{(Add\arabic{longlist})}
\renewcommand\labellonglist{\thelonglist}
\begin{longlist}
\setcounter{longlist}{1}
\item\label{add2}
The bandwidth $h$ is such that (i) $T h^5 \rightarrow\infty$, (ii)
$\frac{\phi_T \log T}{T^{\theta} h^2} = o(1)$ with $\phi_T = \log
\log T$ and $\theta= \min\{ \frac{\beta-4}{\beta}, \frac{\beta
(1-{2}/{s}) - {2}/{s} - 3}{\beta+1} \}$ and (iii) $(T^r
h)^{-1} = o(h^2)$ and $T^{-{r}/({r+1})} = o(h^2)$ with $r = \min\{
\rho,1\}$ and $\rho$ given in~\ref{C1}.
\end{longlist}}

\noindent Condition (ii) is already known from Section
\ref{section-estimation}. As will be seen in Appendix~\ref{appC}, (iii)
ensures that the additional nonstationarity bias is of smaller order
than $O(h^2)$ and can thus be asymptotically neglected. The expressions
for $\beta$ and $\theta$ in~\ref{add1} and~\ref{add2} are calculated as
follows: using the formulas (\ref{beta}) and (\ref{theta}) from
Theorem~\ref{theo-stochastic-part}, we get a pair of expressions for
$\beta$ and $\theta$ for each of the kernel estimates occurring in
(\ref{int-equ}). Combining these expressions yields the formulas in
\ref{add1} and~\ref{add2}.

Under the above assumptions, we can establish the following results,
the proofs of which are given in Appendix~\ref{appC}.
First, the backfitting estimates uniformly converge to the true
component functions at the two-dimensional rates no matter how large
the dimension $d$ of the full regression function.

%th5.1 #&#
%
\begin{theorem}\label{theorem-sbf-1}
Let $I_h = [2C_1h, 1 - 2C_1h]$. Then under~\ref{add1} and
\ref{add2},
%
%e21 #&#
%
\begin{equation}
\label{rates1} \sup_{u,x^j \in I_h} \bigl| \tilde{m}_j
\bigl(u,x^j\bigr) - m_j\bigl(u,x^j\bigr) \bigr| =
O_p \biggl( \sqrt{\frac{\log T}{T h^2}} + h^2 \biggr).
\end{equation}
\end{theorem}
Second, the estimates are asymptotically normal if rescaled appropriately.

%th5.2 #&#
%
\begin{theorem}\label{theorem-sbf-2}
Suppose that~\ref{add1} and~\ref{add2} hold. In addition,
let $\theta> \frac{1}{3}$ and $r > \frac{1}{2}$ to ensure that the
bandwidth $h$ can be chosen to satisfy $\Obs h^6 \rightarrow c_h$ for a
constant $c_h$. Then for any $u,x^1,\ldots, x^d \in(0,1)$,
%
%e22 #&#
%
\begin{equation}
\sqrt{\Obs h^2} \lleft[ %
\matrix{
\tilde{m}_1\bigl(u,x^1\bigr) - m_1
\bigl(u,x^1\bigr)
\cr
\vdots
\cr
\tilde{m}_d
\bigl(u,x^d\bigr) - m_d\bigl(u,x^d\bigr)}
\rright] \convd\normal(B_{u,x},V_{u,x}).
\end{equation}
Here, $V_{u,x}$ is a diagonal matrix whose diagonal entries are given
by the expressions $v_j(u,x^j) = \kappa_0^2 \sigma_j^2(u,x^j) /
p_j(u,x^j)$ with $\kappa_0 = \int K^2(\varphi)\,d\varphi$. Moreover,
the bias term has the form $B_{u,x} = \sqrt{c_h} [\beta_1(u,x^1) -
\gamma_1(u),\ldots, \beta_d(u,x^d) - \gamma_d(u)]^T$. The functions
$\beta_j(u,\cdot)$ in this expression are defined as the minimizers
of the problem
\[
\int\bigl[ \beta(u,x) - b_0 - b_1\bigl(x^1
\bigr) - \cdots- b_d\bigl(x^d\bigr) \bigr]^2
p(u,x) \,dx,
\]
where the minimization runs over all additive functions $b(x) = b_0 +
b_1(x^1) + \cdots+ b_d(x^d)$ with $\int b_j(x^j) p_j(u,x^j) \,dx^j = 0$,
and the function $\beta$ is given in Lem\-ma~\ref{lemma-add-nw-bias} of
Appendix~\ref{appC}. Moreover, the terms $\gamma_j$ can be
characterized by the
equation $\int\alpha_{T,j}(u,x^j)\hat{p}_j(u,x^j) \,dx^j = h^2 \gamma
_j(u) + o_p(h^2)$, where the functions $\alpha_{T,j}$ are again
defined in Lemma~\ref{lemma-add-nw-bias}.
\end{theorem}

%s6 #&#
\section{Concluding remarks}

In this paper, we have studied nonparametric models with a time-varying
regression function and locally stationary covariates. We have
developed a complete asymptotic theory for kernel estimates in these
models. In addition, we have shown that the main assumptions of the
theory are satisfied for a large class of nonlinear autoregressive
processes with a time-varying regression function.

Our analysis can be extended in several directions. An important issue
is bandwidth selection in our framework. As shown in Theorem \ref
{theo-normality}, the asymptotic bias and variance expressions of our
NW estimate are very similar in structure to those from a standard
stationary random design. We thus conjecture that the techniques to
choose the bandwidth in such a design can be adapted to our setting. In
particular, using the formulas for the asymptotic bias and variance
from Theorem~\ref{theo-normality}, it should be possible to select the
bandwidth via plug-in methods.

Another issue concerns forecasting. The convergence results of Theorems~\ref{theo-nw}
and~\ref{theorem-sbf-1} are only valid for rescaled
time lying in a subset $[Ch,1-Ch]$ of the unit interval. For
forecasting purposes, it would be important to provide convergence
rates also in the boundary region $(1-Ch,1]$. This can be achieved by
using boundary-corrected kernels. Another possibility is to work with
one-sided kernels. In both cases, we have to ensure that the kernels
have compact support and are Lipschitz to get the theory to work.

%sA #&#
%
\begin{appendix}\label{app}
%sB #&#
\section{}\label{appA}

In this Appendix, we prove the results on the tvNAR process from
Section~\ref{section-tvNAR}. To shorten notation, we frequently make
use of the abbreviations $\Xvec_{t,T} = \X_{t,T}^{t-d+1}$, $\Xvec_t(u)
= \X_t^{t-d+1}(u)$ and $\epsvec_t = \varepsilon_t^{t-d+1}$.
Moreover, throughout the Appendices, the symbol $C$ denotes a universal
real constant which may take a different value on each occurrence.

\subsection*{Preliminaries}

Before we come to the proofs of the theorems, we state some useful
facts needed for the arguments later on.\vspace*{6pt}

\textsc{Linearization of $m$ and $\sigma$.}
Consider the function $m$. The mean value theorem allows us to write
%
%eB.1 #&#
%
\begin{equation}
\label{mean-val-theo3}\qquad
m \bigl(v,\Xvec_{t-1}(v) \bigr) - m \bigl(u,
\Xvec_{t-1}(u) \bigr) = \Delta_{t,0}^m + \sum
_{j=1}^d \Delta_{t,j}^m
\bigl( \X_{t-j}(v) - \X_{t-j}(u) \bigr),
\end{equation}
where we have used the shorthands $\Delta_{t,0}^m = m(v,\Xvec_{t-1}(v))
- m(u,\Xvec_{t-1}(v))$ and $\Delta_{t,j}^m = \Delta_j^m(u,\Xvec
_{t-1}(u), \Xvec_{t-1}(v))$ for $j=1,\ldots,d$ with the
functions\break $\Delta_j^m (u,y, y') = \int_0^1 \partial_j m(u,y+
s(y' - y)) \,ds$.

The terms $\Delta_{t,j}^m$ have the property that
%
%eB.2 #&#
%
\begin{equation}
\label{hj-bound1} \bigl| \Delta_{t,j}^m \bigr| \le
\Delta_t:= \Delta I\bigl(\|\epsvec_{t-1}\|_{\infty} \le
K_2\bigr) + \delta I\bigl(\|\epsvec_{t-1}\|_{\infty} >
K_2\bigr)
\end{equation}
for $j = 1,\ldots,d$ with $K_2 = (K_1 + \um) / \lsigma$ and $\Delta
\ge\sup_{u,y} |\partial_j m(u,y)|$. This is a straightforward
consequence of the boundedness assumptions on $m$ and $\sigma$. See
the supplement~\cite{Vogt2012} for details.

Repeating the above considerations for the function $\sigma$, we
obtain analogous terms $\Delta_{t,j}^{\sigma}$ that are again bounded
by $\Delta_t$ for $j=1,\ldots,d$.\vspace*{6pt}

\textsc{Recursive formulas for $\X_{t,T}$.}
For the proof of Theorem~\ref{theo-mixing}, we rewrite $\X_{t,T}$ in
a recursive fashion: letting $y_{t-k_1}^{t-k_2}$ and
$e_{t-k_1}^{t-k_2}$ be values of $\X_{t-k_1}^{t-k_2}$ and $\varepsilon
_{t-k_1}^{t-k_2}$, respectively, we recursively define the functions
$m_{t,T}^{(i)}$ by $m_{t,T}^{(0)}(y_{t-1}^{t-d}) = m(\frac{t}{T},
y_{t-1}^{t-d})$
and for $ i \ge1$ by
\begin{eqnarray*}
&&
m_{t,T}^{(i)} \bigl( e_{t-1}^{t-i},y_{t-i-1}^{t-i-d}
\bigr) \\
&&\qquad = m_{t,T}^{(i-1)} \bigl( e_{t-1}^{t-i+1},m_{t-i,T}^{(0)}
\bigl(y_{t-i-1}^{t-i-d}\bigr) + \sigma_{t-i,T}^{(0)}\bigl(y_{t-i-1}^{t-i-d}
\bigr) e_{t-i}, y_{t-i-1}^{t-i-d+1} \bigr).
\end{eqnarray*}
Using analogous recursions for the function $\sigma$, we can
additionally define functions $\sigma_{t,T}^{(i)}$ for $i \ge0$. With
this notation at hand, $\X_{t,T}$ can represented as
\[
\X_{t,T} = m_{t,T}^{(i)} \bigl( \varepsilon_{t-1}^{t-i},
\X_{t-i-1,T}^{t-i-d} \bigr) + \sigma_{t,T}^{(i)}
\bigl( \varepsilon_{t-1}^{t-i},\X_{t-i-1,T}^{t-i-d}
\bigr) \varepsilon_t.
\]
Moreover, for $i \ge d$ we can write
\begin{eqnarray*}
&& m_{t,T}^{(i)} \bigl( e_{t-1}^{t-i},y_{t-i-1}^{t-i-d}
\bigr)
\\
&&\qquad= m \biggl( \frac{t}{T}, m_{t-1,T}^{(i-1)}
\bigl(e_{t-2}^{t-i},y_{t-i-1}^{t-i-d}\bigr) +
\sigma_{t-1,T}^{(i-1)}\bigl(e_{t-2}^{t-i},y_{t-i-1}^{t-i-d}
\bigr) e_{t-1}, \ldots,
\\
&&\qquad\hspace*{40pt} m_{t-d,T}^{(i-d)}\bigl(e_{t-d-1}^{t-i},y_{t-i-1}^{t-i-d}
\bigr) + \sigma_{t-d,T}^{(i-d)}\bigl(e_{t-d-1}^{t-i},y_{t-i-1}^{t-i-d}
\bigr) e_{t-d} \biggr).
\end{eqnarray*}
The term $\sigma_{t,T}^{(i)}( e_{t-1}^{t-i},y_{t-i-1}^{t-i-d} )$ can
be reformulated in the same way.\vspace*{6pt}\vadjust{\goodbreak}

\textsc{Formulas for conditional densities.}
Throughout the Appendix, the symbol $f_{V|W}$ is used to denote the
density of $V$ conditional on $W$. If the residuals $\varepsilon_t$
have a density $f_{\varepsilon}$, then it can be shown that for $1 \le
r \le d$,
%
%eB.3 #&#
%
\begin{equation}
\label{cond-dens} f_{\X_{t,T}|\X_{t-1,T}^{t-r+1},\varepsilon
_{t-r}^{-s},\X
_{-s-1,T}^{-s-d}} \bigl(y_t|y_{t-1}^{t-r+1},e_{t-r}^{-s},z
\bigr) = \frac
{1}{\sigma_{t,T}} f_{\varepsilon} \biggl( \frac{y_t - m_{t,T}}{\sigma
_{t,T}} \biggr).
\end{equation}
Here, $y_t$, $y_{t-1}^{t-r+1}$, $e_{t-r}^{-s}$ and $z$ are values of
$\X_{t,T}$, $\X_{t-1,T}^{t-r+1}$, $\varepsilon_{t-r}^{-s}$ and $\X
_{-s-1,T}^{-s-d}$, respectively. Moreover,
\begin{eqnarray*}
m_{t,T} &=& m \biggl( \frac
{t}{T},y_{t-1}^{t-r+1},m_{t-r,T}^{(t-r+s)}
\bigl(e_{t-r-1}^{-s},z\bigr) + \sigma_{t-r,T}^{(t-r+s)}
\bigl(e_{t-r-1}^{-s},z\bigr) e_{t-r},\ldots,
\\
&&\hspace*{85pt}m_{t-d,T}^{(t-d+s)}\bigl(e_{t-d-1}^{-s},z
\bigr) + \sigma_{t-d,T}^{(t-d+s)}\bigl(e_{t-d-1}^{-s},z
\bigr) e_{t-d} \biggr),
\end{eqnarray*}
and $\sigma_{t,T}$ is defined analogously.

\subsection*{Proof of Theorem \protect\ref{theo-stat}}

Property~\ref{theo-stat-1} follows by standard arguments to be found, for
example, in Chen and Chen~\cite{Chen2000}. Property~\ref{theo-stat-2}
immediately follows with the help of (\ref{cond-dens}). Recalling that
$\X_{t-1,T}^{t-d} = \X_{t-1}^{t-d}(0)$ for $t \le1$,
\ref{theo-stat-3} can again be shown by using~(\ref{cond-dens}).

\subsection*{Proof of Theorem \protect\ref{theo-loc-stat}}

We apply the triangle inequality to get
\[
\bigl| \X_{t,T} - \X_t(u) \bigr| \le\biggl| \X_{t,T} -
\X_t \biggl(\frac{t}{T} \biggr) \biggr| + \biggl| \X_t \biggl(
\frac{t}{T} \biggr) - \X_t(u) \biggr|
\]
and bound the terms $|\X_{t,T} - \X_t(\frac{t}{T})|$ and $|\X_t(\frac
{t}{T}) - \X_t(u)|$ separately. In what follows, we restrict
attention to the term $|\X_t(\frac{t}{T}) - \X_t(u)|$, the arguments
for $|\X_{t,T} - \X_t(\frac{t}{T})|$ being analogous.\vspace*{6pt}

\textsc{Notation.}
Throughout the proof, the symbol $\| z \|$ denotes the Euclidean norm
for vectors $z \in\reals^d$, and $\| A \|$ is the spectral norm for
$d \times d$ matrices $A = (a_{ik})_{i,k = 1,\ldots,d}$. In addition,
$\| A \|_1 = \max_{k=1,\ldots,d} \sum_{j=1}^d |a_{jk}|$.
Furthermore, for $z \in\reals$, we define the family of matrices
\[
B(z) = \pmatrix{ z & \cdots& z & z
\cr
1 & & 0 & 0
\cr
& \ddots& & \vdots
\cr
0 &
& 1 & 0 }.
\]
Finally,\vspace*{2pt} as already noted at the beginning of the Appendix, we make use
of the shorthands $\Xvec_{t,T} = \X_{t,T}^{t-d+1}$, $\Xvec_t(u) = \X
_t^{t-d+1}(u)$ and $\epsvec_t = \varepsilon_t^{t-d+1}$.\vspace*{6pt}

\textsc{Backward iteration.}
By the smoothness conditions on $m$ and $\sigma$,
\[
\X_t \biggl(\frac{t}{T} \biggr) - \X_t(u) =
\bigl( \Delta_{t,0}^m + \Delta_{t,0}^{\sigma}
\varepsilon_t \bigr) + \sum_{j=1}^d
\bigl( \Delta_{t,j}^m + \Delta_{t,j}^{\sigma}
\varepsilon_t \bigr) \biggl( \X_{t-j} \biggl(
\frac{t}{T} \biggr) - \X_{t-j}(u) \biggr)
\]
with $\Delta_{t,0}^m = m(\frac{t}{T},\Xvec_{t-1}(\frac{t}{T})) -
m(u,\Xvec_{t-1}(\frac{t}{T}))$ and $\Delta_{t,j}^m =
\Delta_j^m(u,\Xvec_{t-1}(u),\break \Xvec_{t-1}(\frac{t}{T}))$ for
$j=1,\ldots,d$ as introduced in (\ref{mean-val-theo3}). The terms
$\Delta_{t,j}^{\sigma}$ for $j = 0,\ldots,d$ are defined analogously.
In matrix notation, we obtain
%
%eB.4 #&#
%
\begin{equation}
\label{iterate1} \Xvec_t \biggl(\frac{t}{T} \biggr) -
\Xvec_t(u) = \Dvec_t \biggl( \Xvec_{t-1}
\biggl(\frac{t}{T} \biggr) - \Xvec_{t-1}(u) \biggr) +
\Evec_t
\end{equation}
with $\Evec_t = (\Delta_{t,0}^m + \Delta_{t,0}^{\sigma} \varepsilon
_t,0,\ldots,0)^T$ and
\[
\Dvec_t = \pmatrix{ \Delta_{t,1}^m +
\Delta_{t,1}^{\sigma} \varepsilon_t & \cdots&
\Delta_{t,d-1}^m + \Delta_{t,d-1}^{\sigma}
\varepsilon_t & \Delta_{t,d}^m +
\Delta_{t,d}^{\sigma} \varepsilon_t
\cr
1 & & 0 & 0
\cr
& \ddots& & \vdots
\cr
0 & & 1 & 0 }.
\]
Iterating (\ref{iterate1}) $n$ times yields
\begin{eqnarray*}
\biggl\| \Xvec_t \biggl(\frac{t}{T} \biggr) - \Xvec_t(u)
\biggr\| & \le&\| \Evec_t \| + \Biggl\| \sum_{r=0}^{n-1}
\prod_{k=0}^r \Dvec_{t-k}
\Evec_{t-r-1} \Biggr\|
\\
&&{} + \Biggl\| \prod_{k=0}^n
\Dvec_{t-k} \biggl( \Xvec_{t-n-1} \biggl(\frac{t}{T}
\biggr) - \Xvec_{t-n-1}(u) \biggr) \Biggr\|.
\end{eqnarray*}
Note that the rescaled time argument $\frac{t}{T}$ plays the same role
as the argument $u$ and thus remains fixed when iterating backward.
Next define matrices $\DDvec_t$ by
%
%eB.5 #&#
%
\begin{equation}
\label{iteration-matrix} \DDvec_t = \bigl(1 + |\varepsilon_t|\bigr)
B(\Delta_t)
\end{equation}
with $\Delta_t = \Delta I(\|\epsvec_{t-1}\|_{\infty} \le K_2) + \delta
I(\|\epsvec_{t-1}\|_{\infty} > K_2)$. As shown in the preliminaries
section of the Appendix, $|\Delta_{t,j}^m + \Delta_{t,j}^{\sigma}
\varepsilon_t| \le\Delta_t (1 + |\varepsilon_t|)$ for $j=1,\ldots,d$.
Therefore, the entries of the matrix $\DDvec_t$ are all weakly larger
in absolute value than those of $\Dvec_t$. This implies that $\| \prod
_{k=0}^n \Dvec_{t-k} z \| \le\| \prod_{k=0}^n \DDvec_{t-k} |z| \|$
with $z = (|z_1|,\ldots,|z_d|)$. Using this together with the
boundedness of $m$ and $\sigma$ and the fact that $|\Delta_{t,0}^m +
\Delta_{t,0}^{\sigma} \varepsilon_t| \le C|\frac{t}{T} - u| (1 +
|\varepsilon_t|)$, we finally arrive at
\[
\Biggl\| \Xvec_t \biggl(\frac{t}{T} \biggr) - \Xvec_t(u)
\Biggr\| \le\biggl| \frac{t}{T} - u \biggr| V_{t,n} + R_{t,n}
\]
with
\begin{eqnarray*}
V_{t,n} & = & C \bigl(1 + |\varepsilon_t|\bigr) + C \sum
_{r=0}^{n-1} \bigl(1 + |\varepsilon_{t-r-1}| \bigr) \Biggl\|
\prod_{k=0}^r \DDvec_{t-k} \Biggr\|,
\\
R_{t,n} & = & C \bigl(1 + \|\epsvec_{t-n-1}\|\bigr) \Biggl\| \prod
_{k=0}^n \DDvec_{t-k} \Biggr\|.
\end{eqnarray*}

\textsc{Bounding $V_{t,n}$ and $R_{t,n}$.}
The convergence behavior of $V_{t,n}$ and $R_{t,n}$ for $n \rightarrow
\infty$ mainly depends on the properties of the product $\| {\prod}
_{k=0}^n \DDvec_{t-k} \|$. The behavior of the latter is
described by the following lemma.
%
%leB.1 #&#
%
\begin{lemma}\label{lemmaA1}
If $\delta$ is sufficiently small, in particular, if it satisfies
(\ref{delta}), then there exists a constant $\rho> 0$ such that for
some $\gamma< 1$,
%
%eB.6 #&#
%
\begin{equation}
\label{lemmaA1-statement} \ex\Biggl[ \Biggl\| \prod_{k=0}^n
\DDvec_{t-k} \Biggr\|^{\rho} \Biggr] \le C \gamma^n.
\end{equation}
\end{lemma}
The proof of Lemma~\ref{lemmaA1} is postponed until the arguments for
Theorem~\ref{theo-loc-stat} are completed. The following statement is
a direct consequence of Lemma~\ref{lemmaA1}.
\begin{longlist}[(R)]
\item[(R)] There exists a constant $\rho> 0$ such that $\ex[ R_{t,n}^{\rho
} ] \le C \gamma^n$ for some $\gamma< 1$. In particular, $R_{t,n}
\convas0$ as $n \rightarrow\infty$.
\end{longlist}
In addition, it holds that:
\begin{longlist}[(V)]
\item[(V)] $V_{t,n} \le V_t$, where the variables $V_t$ have the property
that $\ex[ V_t^{\rho} ] \le C$ for a positive constant $\rho< 1$ and
all $t$.
\end{longlist}
This can be seen as follows. First note that
\[
V_{t,n} \le C \bigl(1 + |\varepsilon_t|\bigr) + \sum
_{r=0}^{n-1} R_{t,r} \le V_t:=
C \bigl(1 + |\varepsilon_t|\bigr) + \sum_{r=0}^{\infty}
R_{t,r}.
\]
Using the monotone convergence theorem and Lo\`eve's inequality with
$\rho< 1$, we obtain $\ex[ V_t^{\rho} ] \le C \ex(1 + |\varepsilon
_t|)^{\rho} + \sum_{r=0}^{\infty} \ex[R_{t,r}^{\rho}]$.
As the right-hand side of the previous inequality is finite by (R), we
arrive at (V).

(R) and (V) imply that $|\X_t(\frac{t}{T}) - \X_t(u)| \le|\frac
{t}{T} - u| V_t$ a.s. with variables $V_t$ whose $\rho$th moment is
uniformly bounded by some finite constant $C$. An analogous result can
be derived for $|\X_{t,T} - \X_t(\frac{t}{T})|$. This completes the
proof.
\begin{pf*}{Proof of Lemma~\ref{lemmaA1}}
We want to show that the $\rho$th moment of the product $\| {\prod}
_{k=0}^n \DDvec_{t-k} \|$ converges exponentially fast to
zero as $n \rightarrow\infty$. This is a highly nontrivial problem,
and as far as we can see, it cannot be solved by simply adapting
techniques from related papers on models with time-varying
coefficients. The problem is that the techniques used therein are
either tailored to products of deterministic matrices (see, e.g.,
Proposition 13 in Moulines et al.~\cite{Moulines2005}) or they
heavily draw on the independence of the random matrices involved (see,
e.g., Proposition~2.1 in Subba Rao~\cite{SubbaRao2006}).

We now describe our proving strategy in detail. To start with, we
replace the spectral norm \mbox{$\| \cdot\|$} in (\ref{lemmaA1-statement})
by the norm \mbox{$\| \cdot\|_1$} which is much easier to handle. As these
two norms are equivalent, there exists a finite constant $C$ such that
$ \| {\prod}_{k=0}^n \DDvec_{t-k} \| \le C \mathcal{B}_n$
with $\mathcal{B}_n = \| {\prod}_{k=0}^n \DDvec_{t-k} \|_1$.
Next, we split up the term $\mathcal{B}_n$ into two parts,
\[
\mathcal{B}_n = I_n \mathcal{B}_n + (1 -
I_n) \mathcal{B}_n =: \mathcal{B}_{n,1} +
\mathcal{B}_{n,2},
\]
where $I_n = I(\sum_{k=0}^n J_k > \kappa n)$ with $J_k = I(
\min_{l=1,\ldots,d} |\varepsilon_{t-k-l}| \le K_2)$ and a constant
$0 < \kappa< 1$ to be specified later on. Lemma~\ref{lemmaA1} is a
direct consequence of the following two facts:
\begin{longlist}[(ii)]
\item[(i)] There exists a constant $\rho> 0$ such that $\ex[
\mathcal{B}_{n,1}^{\rho} ] \le C \gamma^n$ for some $\gamma< 1$.
\item[(ii)] $\ex[\mathcal{B}_{n,2}] \le C \gamma^n$ for some
$\gamma< 1$.
\end{longlist}

We start with the proof of (i). Letting $\phi_n = \lambda^n$ with
some positive constant $\lambda< 1$, we can write
\begin{eqnarray*}
\ex\bigl[ \mathcal{B}_{n,1}^{\rho} \bigr] & = & \ex\bigl[ I(
\mathcal{B}_{n,1} > \phi_n ) \mathcal{B}_{n,1}^{\rho}
\bigr] + \ex\bigl[ I( \mathcal{B}_{n,1} \le\phi_n )
\mathcal{B}_{n,1}^{\rho} \bigr]
\\
& \le&\bigl( \ex\bigl[ \mathcal{B}_{n,1}^{2\rho} \bigr] \pr(
\mathcal{B}_{n,1} > \phi_n ) \bigr)^{1/2} +
\phi_n^{\rho}.
\end{eqnarray*}
It is easy to see that $\ex[ \mathcal{B}_{n,1}^{2\rho} ] \le C^{\rho
n}$ for a sufficiently large constant $C$, where $C^{\rho}$ can be
made arbitrarily close to one by choosing $\rho> 0$ small enough. To
show (i), it thus suffices to verify that
%
%eB.7 #&#
%
\begin{equation}
\label{Bn1a} \pr( \mathcal{B}_{n,1} > \phi_n ) \le C
\gamma^n \qquad\mbox{for some } \gamma< 1.
\end{equation}
For the proof of (\ref{Bn1a}), we write
\[
\pr( \mathcal{B}_{n,1} > \phi_n ) \le\pr(
I_n > 0 ) = \pr\Biggl( \sum_{k=0}^n
\bigl(J_k - \ex[J_k]\bigr) > \kappa_0 n
\Biggr)
\]
with $\kappa_0:= \kappa- \ex[J_k]$. As the variables $\varepsilon_t$
have an everywhere positive density by assumption, the expectation
$\ex[J_k]$ is strictly smaller than one. We can thus choose $0 <
\kappa< 1$ slightly larger than $\ex[J_k]$ to get that $0 < \kappa_0
< 1$. As the variables $J_k - \ex[J_k]$ for $k=0,\ldots,n$ are
$2d$-dependent, a simple blocking argument together with Hoeffding's
inequality shows that
\[
\pr\Biggl( \sum_{k=0}^{n}
\bigl(J_k - \ex[J_k]\bigr) > \kappa_0 n
\Biggr) \le C \gamma^n
\]
for some $\gamma< 1$. This yields (\ref{Bn1a}) and thus completes the
proof of (i).

Let us now turn to the proof of (ii). We have that
\[
\mathcal{B}_{n,2} = (1 - I_n) \prod
_{k=0}^n \bigl(1 - |\varepsilon_{t-k}|\bigr) \Biggl\|
\prod_{k=0}^n \DDvec(\Delta_{t-k})
\Biggr\|_1.
\]
The random matrix $B(\Delta_{t-k})$ in the above expression can only
take two forms: if $\| \underline{\varepsilon}_{t-k-1} \|_{\infty} >
K_2$, it equals $B(\delta)$, and if $\| \underline{\varepsilon
}_{t-k-1} \|_{\infty} \le K_2$, it equals $B(\Delta)$. Moreover, if
$\min_{l=1,\ldots,d} |\varepsilon_{t-k-l}| > K_2$, it holds that $\|
\epsvec_{t-k-l} \|_{\infty} > K_2$ for all $l=1,\ldots,d$ and thus
$\prod_{l=0}^{d-1} B(\Delta_{t-k-l}) = B(\delta)^d$.
Importantly, the term $\mathcal{B}_{n,2}$ is unequal to zero only if
$I_n = 0$, that is, only if $\min_{l=1,\ldots,d} |\varepsilon_{t-k-l}|
> K_2$ for at least $(1 - \kappa) n$ terms. From this, we
can infer that
%
%eB.8 #&#
%
\begin{equation}
\label{mean1} \ex[ \mathcal{B}_{n,2} ] \le\ex\Biggl[ \prod
_{k=0}^n \bigl(1 + |\varepsilon_{t-k}|\bigr) \Biggr]
\bigl\| B(\Delta) \bigr\|_1^{\kappa n} \bigl\| B(\delta)^d
\bigr\|_1^{{(1 - \kappa) n}/{d}}.
\end{equation}
By direct calculations, we can verify that $\| B(\delta)^d \|_1 \le
C_d \delta$ with the constant $C_d = \sum_{l=0}^{d-1} \sum_{k=0}^l
{l\choose k}$ that only depends on the dimension
$d$. Moreover, $\|B(\Delta)\|_1 \le(\Delta+ 1)$. Plugging this into
(\ref{mean1}) yields
\[
\ex[ \mathcal{B}_{n,2} ] \le\bigl(1 + \ex|\varepsilon_0|\bigr)
\bigl[ \bigl(1 + \ex|\varepsilon_0|\bigr) (\Delta+ 1)^{\kappa}
(C_d \delta)^{{(1 - \kappa)}/{d}} \bigr]^n.
\]
Straightforward calculations show that the term in square brackets is
strictly smaller than one for
%
%eB.9 #&#
%
\begin{equation}
\label{delta} \delta< \bigl[ \bigl(1 + \ex|\varepsilon_0|\bigr)^{{d}/({1-\kappa})}
(\Delta+ 1)^{{\kappa d}/({1 - \kappa})} C_d \bigr]^{-1}.
\end{equation}
Assuming that $\delta$ satisfies the above condition, we thus arrive
at (ii).
\end{pf*}

\subsection*{Proof of Theorem \protect\ref{theo-smooth-dens}}

The proof can be found in the supplement~\cite{Vogt2012}.

\subsection*{Proof of Theorem \protect\ref{theo-mixing}}

To start with, note that the process $\{\X_{t,T}\}$ is $d$-Markovian.
This implies that
\[
\beta(k) = \sup_{T \in\integers} \sup_{t \in\integers} \beta\bigl
(\sigma(
\Xvec_{t-k,T}),\sigma(\Xvec_{t+d-1,T}) \bigr)
\]
with
\[
\beta\bigl(\sigma(\Xvec_{t-k,T}),\sigma(\Xvec_{t+d-1,T}) \bigr) =
\ex\Bigl[ \sup_{\sset\in\sigma(\Xvec_{t+d-1,T})} \bigl| \pr(\sset) - \pr
\bigl(\sset|\sigma(
\Xvec_{t-k,T})\bigr) \bigr| \Bigr].
\]
In the following, we bound the expression $|\pr(\sset) - \pr(\sset
|\sigma(\Xvec_{t-k,T}))|$ for arbitrary sets $\sset\in\sigma(\Xvec
_{t+d-1,T})$. This provides us with a bound for the mixing coefficients
$\beta(k)$ of the process $\{\X_{t,T}\}$.

We use the following notation: throughout the proof, we let $y=
y_{t+d-1}^t$, $e = e_{t-1}^{t-k+1}$ and $z = z_{t-k}^{t-k-d+1}$ be
values of $\Xvec_{t+d-1,T}$, $\varepsilon_{t-1}^{t-k+1}$ and
$\Xvec_{t-k,T}$, respectively. Moreover, we use the shorthand
\[
f_j(y_{t+j}|z) = f_{\X_{t+j,T}|\X_{t+j-1,T}^t,\varepsilon
_{t-1}^{t-k+1},\Xvec_{t-k,T}}\bigl(y_{t+j}|y_{t+j-1}^t,e,z
\bigr)
\]
for $j=0,\ldots,d-1$, where we suppress the dependence on the
arguments $y_{t+j-1}^t$ and $e$ in the notation.
Finally, note that by (\ref{cond-dens}), the above conditional density
can be expressed in terms of the error density $f_{\varepsilon}$ as
%
%eB.10 #&#
%
\begin{equation}
\label{cond-dens'} f_j(y_{t+j}|z) =
\frac{1}{\sigma_{t,T,j}(z)} f_{\varepsilon} \biggl(\frac{y_{t+j} -
m_{t,T,j}(z)}{\sigma_{t,T,j}(z)} \biggr)
\end{equation}
with
\begin{eqnarray*}
m_{t,T,j}(z) & = & m \biggl(\frac
{t+j}{T},y_{t+j-1}^t,m_{t-1,T}^{(k-2)}
\bigl(e_{t-2}^{t-k+1},z\bigr) + \sigma_{t-1,T}^{(k-2)}
\bigl(e_{t-2}^{t-k+1},z\bigr) e_{t-1},\ldots,
\\
&&\qquad\hspace*{5.2pt} m_{t+j-d,T}^{(k-j+d-1)}\bigl(e_{t+j-d-1}^{t-k+1},z
\bigr) + \sigma_{t+j-d,T}^{(k-j+d-1)}\bigl(e_{t+j-d-1}^{t-k+1},z
\bigr) e_{t+j-d} \biggr)
\end{eqnarray*}
and $\sigma_{t,T,j}(z)$ defined analogously. The functions
$m_{t-1,T}^{(k-2)}$, $\sigma_{t-1,T}^{(k-2)}, \ldots$ were introduced
in the preliminaries section of the Appendix.

With this notation at hand, we can write
\begin{eqnarray*}
&& \pr\bigl(\sset|\sigma(\Xvec_{t-k,T})\bigr)
\\
&&\qquad= \ex\bigl[ \ex\bigl[ I(\Xvec_{t+d-1,T} \in\sset) |
\varepsilon_{t-1}^{t-k+1},\Xvec_{t-k,T} \bigr] |
\Xvec_{t-k,T} \bigr]
\\
&&\qquad= \int I(y\in\sset) f_{\Xvec_{t+d-1,T}|\varepsilon
_{t-1}^{t-k+1},\Xvec_{t-k,T}}(y|e,\Xvec_{t-k,T}) \prod
_{l=1}^{k-1} f_{\varepsilon}(e_{t-l}) \,de \,dy
\\
&&\qquad= \int I(y\in\sset) \prod_{j=0}^{d-1}
f_j(y_{t+j}|\Xvec_{t-k,T}) \prod
_{l=1}^{k-1} f_{\varepsilon}(e_{t-l}) \,de \,dy
\end{eqnarray*}
and likewise
\[
\pr(\sset) = \int I(y\in\sset) \prod_{j=0}^{d-1}
f_j(y_{t+j}|z) \prod_{l=1}^{k-1} f_{\varepsilon
}(e_{t-l}) f_{\Xvec_{t-k,T}}(z) \,de \,dz \,dy.
\]
Using the shorthand $\underline{\X} = \Xvec_{t-k,T}$, we thus arrive at
\begin{eqnarray*}
&& \bigl| \pr(\sset) - \pr\bigl(\sset|\sigma(\underline{\X})\bigr) \bigr|
\\
&&\qquad\le\int\underbrace{ \Biggl[ \int\Biggl| \prod_{j=0}^{d-1}
f_j(y_{t+j}|z) - \prod_{j=0}^{d-1}
f_j(y_{t+j}|\underline{\X}) \Biggr| \,dy \Biggr]}_{=: (*)}
\prod_{l=1}^{k-1} f_{\varepsilon
}(e_{t-l})
f_{\underline{\X}}(z) \,de \,dz.
\end{eqnarray*}
We next consider $(*)$ more closely. A telescoping argument together
with Fubini's theorem yields that
\begin{eqnarray*}
(*) & \le&\sum_{i=0}^{d-1} \int\Biggl[
\prod_{j=0}^{i-1} f_j(y_{t+j}|
\underline{\X}) \bigl| f_i(y_{t+i}|z) - f_i(y_{t+i}|
\underline{\X}) \bigr| \prod_{j=i+1}^{d-1}
f_j(y_{t+j}|z) \Biggr] \,dy
\\
& = &\sum_{i=0}^{d-1} \int\Biggl[ \int
\Biggl[ \int\prod_{j=i+1}^{d-1}
f_j(y_{t+j}|z) \,dy_{t+d-1} \cdots dy_{t+i+1}
\Biggr]
\\
&&\hspace*{79pt}{} \times\bigl| f_i(y_{t+i}|z) - f_i(y_{t+i}|
\underline{\X}) \bigr| \,dy_{t+i} \Biggr] \\
&&\hspace*{25.5pt}{}\times\prod_{j=0}^{i-1}
f_j(y_{t+j}|\underline{\X}) \,dy_{t+i-1} \cdots
dy_t
\\
& \le& \sum_{i=0}^{d-1} \int\underbrace{
\biggl[ \int\bigl| f_i(y_{t+i}|z) - f_i(y_{t+i}|
\underline{\X}) \bigr| \,dy_{t+i} \biggr]}_{=:(**)} \prod
_{j=0}^{i-1} f_j(y_{t+j}|
\underline{\X}) \,dy_{t+i-1} \cdots dy_t,
\end{eqnarray*}
where the last inequality exploits the fact that
\[
\int\prod
_{j=i+1}^{d-1} f_j(y_{t+j}|z) \,dy_{t+d-1} \cdots dy_{t+i+1}
\]
is a conditional probability and thus almost surely bounded by one.
Using formula (\ref{cond-dens'}) together with~\ref{E3}, it is
straightforward to see that
\begin{eqnarray*}
(**) & = & \int\biggl| \frac{1}{\sigma_{t,T,i}(z)} f_{\varepsilon} \biggl
(\frac
{y_{t+i} - m_{t,T,i}(z)}{\sigma_{t,T,i}(z)}
\biggr) \\
&&\hspace*{11.6pt}{} - \frac
{1}{\sigma_{t,T,i}(\underline{\X})} f_{\varepsilon} \biggl(\frac
{y_{t+i} - m_{t,T,i}(\underline{\X})}{\sigma_{t,T,i}(\underline{\X
})} \biggr) \biggr|
\,dy_{t+i}
\\
& \le& C \bigl( \bigl| m_{t,T,i}(z) - m_{t,T,i}(\underline{\X}) \bigr| + \bigl|
\sigma_{t,T,i}(z) - \sigma_{t,T,i}(\underline{\X}) \bigr| \bigr)
\\
& \le& C (2 \um+ 2 \usigma) \bigl( \bigl| m_{t,T,i}(z) - m_{t,T,i}(
\underline{\X}) \bigr| + \bigl| \sigma_{t,T,i}(z) - \sigma_{t,T,i}(
\underline{\X}) \bigr| \bigr)^p,
\end{eqnarray*}
where $p$ is some constant with $0 < p < 1$. Iterating backward $n \le
k-2d$ times in the same way as in Theorem~\ref{theo-loc-stat}, we can
further show that
%
%eB.11 #&#
%
\begin{eqnarray}\label{mixing-star3}
&&
\bigl| m_{t,T,i}(z) - m_{t,T,i}(\underline{\X}) \bigr| + \bigl|
\sigma_{t,T,i}(z) - \sigma_{t,T,i}(\underline{\X}) \bigr|
\nonumber\\[-8pt]\\[-8pt]
&&\qquad\le C \sum_{j=1}^{d-i} \Biggl\| \prod
_{m=0}^n \DDvec_{t-j-m} \Biggr\| \bigl( 1 +
\bigl\|e_{t-j-n-1}^{t-j-n-d}\bigr\| \bigr),
\nonumber
\end{eqnarray}
where $\| \cdot\|$ denotes the Euclidean norm for vectors and the
spectral norm for matrices. The matrix $\DDvec_t$ was introduced in
(\ref{iteration-matrix}). Note that $\DDvec_t$ was defined there in
terms of the random vector $\varepsilon_t^{t-d}$. Slightly abusing
notation, we here use the symbol $\DDvec_t$ to denote the matrix with
$\varepsilon_t^{t-d}$ replaced by the realization $e_t^{t-d}$. Keeping
in mind that the matrix $\DDvec_t$ only depends on the residual values
$e_t^{t-d}$, we can plug (\ref{mixing-star3}) into the bound for
$(**)$ and insert this into the bound for $(*)$ to arrive at
\[
(*) \le C \Biggl( \sum_{j=1}^d \Biggl\| \prod
_{m=0}^n \DDvec_{t-j-m} \Biggr\| \bigl(
1 + \bigl\|e_{t-j-n-1}^{t-j-n-d}\bigr\| \bigr) \Biggr)^p.
\]
As a consequence,
\[
\bigl| \pr(\sset) - \pr\bigl(\sset|\sigma(\underline{\X})\bigr) \bigr|
\le C \ex
\Biggl(
\sum_{j=1}^d \Biggl\| \prod
_{m=0}^n \DDvec_{t-j-m} \Biggr\| \bigl( 1 + \bigl\|
\varepsilon_{t-j-n-1}^{t-j-n-d}\bigr\| \bigr) \Biggr)^p.
\]
Using the arguments from Lemma~\ref{lemmaA1}, we can show that for $p
> 0$ sufficiently small, the expectation on the right-hand side is
bounded by $C \lambda^n$ for some positive constant $\lambda< 1$.
Choosing $n = k-2d$, for instance, we thus arrive at
\[
\bigl| \pr(\sset) - \pr\bigl(\sset|\sigma(\Xvec_{t-k,T})\bigr) \bigr| \le C
\lambda^{k-(d+1)} \le C \gamma^k
\]
for some constant $\gamma< 1$. This immediately implies that $\beta
(k) \le C \gamma^k$.

%sC #&#
\section{}\label{appB}
In this Appendix, we prove the results of Section \ref
{section-estimation}. Before we turn to the proofs, we state two
auxiliary lemmas which are repeatedly used throughout the Appendix. The
proofs are straightforward and thus omitted.
%
%leC.1 #&#
%
\begin{lemma}\label{lemmaB1}
Suppose the kernel $K$ satisfies~\ref{C6} and let $I_h = [C_1 h,1
- C_1 h]$. Then for $k=0,1,2$,
\begin{eqnarray*}
&&\sup_{u \in I_h} \Biggl| \frac{1}{Th} \sum_{t=1}^T
K_h \biggl(u - \frac{t}{T} \biggr) \biggl( \frac{u - {t}/{T}}{h}
\biggr)^k - \int_0^1
\frac{1}{h} K_h(u - \varphi) \biggl( \frac{u - \varphi
}{h}
\biggr)^k \,d\varphi\Biggr| \\
&&\qquad= O \biggl(\frac{1}{Th^2} \biggr).
\end{eqnarray*}
\end{lemma}
%
%leC.2 #&#
%
\begin{lemma} \label{lemmaB4}
Suppose $K$ satisfies~\ref{C6} and let $g\dvtx[0,1] \times\reals^d
\rightarrow\reals$, $(u,x) \mapsto g(u,x)$ be continuously
differentiable w.r.t. $u$. Then for any compact set $S \subset\reals^d$,
\[
\sup_{u \in I_h, x \in S} \Biggl| \frac{1}{T h} \sum_{t=1}^T
K_h \biggl(u - \frac{t}{T} \biggr) g \biggl(
\frac{t}{T},x \biggr) - g(u,x) \Biggr| = O \biggl(\frac{1}{Th^2} \biggr) +
o(h).
\]
\end{lemma}

\subsection*{Proof of Theorem \protect\ref{theo-stochastic-part}}

To show the result, we use a blocking argument together with an
exponential inequality for mixing arrays, thus following the common
proving strategy to be found, for example, in Bosq~\cite{Bosq1998},
Masry~\cite{Masry1996} or Hansen~\cite{Hansen2008}. In particular, we
go along the lines of Hansen's proof of Theorem~2 in~\cite
{Hansen2008}, modifying his arguments to allow for local
stationarity
in the data. A detailed version of the arguments can be found in the
supplement~\cite{Vogt2012}.

\subsection*{Proof of Theorem \protect\ref{theo-nw}}

We write
\[
\hat{m}(u,x) - m(u,x) = \frac{1}{\hat{f}(u,x)} \bigl( \hat{g}^V(u,x) +
\hat{g}^B(u,x) - m(u,x) \hat{f}(u,x) \bigr)\vadjust{\goodbreak}
\]
with
\begin{eqnarray*}
\hat{f}(u,x) & = & \frac{1}{T h^{d+1}} \sum_{t=1}^T
K_h \biggl(u - \frac{t}{T} \biggr) \prod
_{j=1}^d K_h\bigl(x^j -
X_{t,T}^j\bigr),
\\
\hat{g}^V(u,x) & = & \frac{1}{T h^{d+1}} \sum
_{t=1}^T K_h \biggl(u -
\frac{t}{T} \biggr) \prod_{j=1}^d
K_h\bigl(x^j - X_{t,T}^j\bigr)
\varepsilon_{t,T},
\\
\hat{g}^B(u,x) & = & \frac{1}{T h^{d+1}} \sum
_{t=1}^T K_h \biggl(u -
\frac{t}{T} \biggr) \prod_{j=1}^d
K_h\bigl(x^j - X_{t,T}^j\bigr) m
\biggl(\frac{t}{T},X_{t,T} \biggr).
\end{eqnarray*}
We first derive some intermediate results for the above expressions:
{\renewcommand\thelonglist{(\roman{longlist})}
\renewcommand\labellonglist{\thelonglist}
\begin{longlist}
\item\label{conv-a} By Theorem~\ref{theo-stochastic-part} with
$W_{t,T} = \varepsilon_{t,T}$,
\[
\sup_{u \in[0,1], x \in S} \bigl| \hat{g}^V(u,x) \bigr| = O_p \biggl(
\sqrt{\frac{\log T}{T h^{d+1}}} \biggr).
\]
\item\label{conv-b} Applying the arguments for Theorem \ref
{theo-stochastic-part} to $\hat{g}^B(u,x) - m(u,x) \hat{f}(u,x)$ yields
\begin{eqnarray*}
&& \sup_{u \in[0,1], x \in S} \bigl| \hat{g}^B(u,x) - m(u,x)
\hat{f}(u,x)\\
&&\qquad\quad\hspace*{12.5pt}{} - \ex\bigl[ \hat{g}^B(u,x) - m(u,x) \hat
{f}(u,x) \bigr] \bigr|\\
&&\qquad =
O_p \biggl( \sqrt{\frac{\log T}{T h^{d+1}}} \biggr).
\end{eqnarray*}
\item\label{conv-c} It holds that
\begin{eqnarray*}
&& \sup_{u \in I_h, x \in S} \bigl| \ex\bigl[ \hat{g}^B(u,x) - m(u,x) \hat
{f}(u,x) \bigr] \bigr|
\\
&&\qquad= h^2 \frac{\kappa_2}{2} \sum_{i=0}^d
\bigl( 2 \,\partial_i m(u,x) \,\partial_i f(u,x) +
\partial_{ii}^2 m(u,x) f(u,x) \bigr)\\
&&\qquad\quad{} + O \biggl(
\frac{1}{T^r h^d} \biggr) + o\bigl(h^2\bigr)
\end{eqnarray*}
with $r = \min\{\rho,1\}$. The proof is postponed until the arguments
for Theorem~\ref{theo-nw} are completed.
\item\label{conv-d} We have that
\[
\sup_{u \in I_h, x \in S} \bigl| \hat{f}(u,x) - f(u,x) \bigr| = o_p(1).
\]
For the proof, we split up the term $\hat{f}(u,x) - f(u,x)$ into a
variance part $\hat{f}(u,x) - \ex\hat{f}(u,x)$ and a bias part $\ex
\hat{f}(u,x) - f(u,x)$. Applying Theorem~\ref{theo-stochastic-part}
with $W_{t,T}=1$ yields that the variance part is $o_p(1)$ uniformly in
$u$. The bias part can be analyzed by a simplified version of the
arguments used to prove~\ref{conv-c}.
\end{longlist}}

\noindent Combining the intermediate results~\ref{conv-a}--\ref{conv-c}, we
arrive at
\begin{eqnarray*}
&& \sup_{u \in I_h, x \in S} \bigl| \hat{m}(u,x) - m(u,x) \bigr|
\\[-2pt]
&&\qquad\le\bigl( \sup\hat{f}(u,x)^{-1} \bigr) \bigl( \sup\bigl| \hat
{g}^V(u,x) \bigr| + \sup\bigl| \hat{g}^B(u,x) - m(u,x)
\hat{f}(u,x) \bigr| \bigr)
\\[-2pt]
&&\qquad= \bigl(\sup\hat{f}(u,x)^{-1} \bigr) O_p \biggl( \sqrt{
\frac
{\log T}{T h^{d+1}}} + \frac{1}{T^r h^d} + h^2 \biggr)
\end{eqnarray*}
with $r = \min\{\rho,1\}$. Moreover,~\ref{conv-d} and the
condition that $\inf_{u \in[0,1], x \in S} f(u,\break x) > 0$ immediately
imply that $\sup\hat{f}(u,x)^{-1} = O_p(1)$. This completes the proof.

\begin{pf*}{Proof of~\ref{conv-c}}
Let $\bar{K}\dvtx\reals\rightarrow\reals$ be a Lipschitz continuous
function with support $[-q C_1, q C_1]$ for some $q > 1$. Assume that
$\bar{K}(x) = 1$ for all $x \in[-C_1,C_1]$ and write $\bar{K}_h(x) =
\bar{K}(\frac{x}{h})$. Then
\[
\ex\bigl[ \hat{g}^B(u,x) - m(u,x) \hat{f}(u,x) \bigr] =
Q_1(u,x) + \cdots+ Q_4(u,x)
\]
with
\[
Q_i(u,x) = \frac{1}{Th^{d+1}} \sum_{t=1}^T
K_h \biggl(u - \frac
{t}{T} \biggr) q_i(u,x)
\]
and
\begin{eqnarray*}
q_1(u,x) & = & \ex\Biggl[ \prod_{j=1}^d
\bar{K}_h\bigl(x^j - X_{t,T}^j
\bigr) \Biggl\{ \prod_{j=1}^d
K_h\bigl(x^j - X_{t,T}^j\bigr)
\\
&&\hspace*{101pt}{} - \prod_{j=1}^d K_h
\biggl(x^j - X_t^j\biggl( {\frac{t}{T}}
\biggr) \biggr) \Biggr\}\\
&&\hspace*{95.4pt}{}\times \biggl\{ m \biggl( {\frac{t}{T}},X_{t,T}
\biggr) - m(u,x) \biggr\} \Biggr],
\\
q_2(u,x) & = & \ex\Biggl[ \prod_{j=1}^d
\bar{K}_h\bigl(x^j - X_{t,T}^j
\bigr) \prod_{j=1}^d K_h
\biggl(x^j - X_t^j\biggl( {\frac{t}{T}}
\biggr) \biggr)
\\
&&\hspace*{34.8pt}{} \times\biggl\{ m \biggl( {\frac
{t}{T}},X_{t,T} \biggr) - m
\biggl( {\frac
{t}{T}},X_t\biggl( {\frac{t}{T}}\biggr)
\biggr) \biggr\} \Biggr],
\\
q_3(u,x) & = & \ex\Biggl[ \Biggl\{ \prod_{j=1}^d
\bar{K}_h\bigl(x^j - X_{t,T}^j
\bigr) - \prod_{j=1}^d
\bar{K}_h \biggl(x^j - X_t^j
\biggl( {\frac{t}{T}}\biggr) \biggr) \Biggr\}
\\[-2pt]
&&\hspace*{11pt}{} \times\prod_{j=1}^d K_h
\biggl(x^j - X_t^j\biggl( {\frac{t}{T}}
\biggr) \biggr) \biggl\{ m \biggl( {\frac{t}{T}},X_t\biggl( {
\frac{t}{T}}\biggr) \biggr) - m(u,x) \biggr\}\Biggr],
\\[-2pt]
q_4(u,x) & = & \ex\Biggl[ \prod_{j=1}^d
K_h \biggl(x^j - X_t^j\biggl( {
\frac{t}{T}}\biggr) \biggr) \biggl\{ m \biggl( {\frac{t}{T}},X_t
\biggl( {\frac{t}{T}}\biggr) \biggr) - m(u,x) \biggr\} \Biggr].
\end{eqnarray*}
We first consider $Q_1(u,x)$. As the kernel $K$ is bounded, we can use
a telescoping argument to get that $| \prod_{j=1}^d K_h(x^j
- X_{t,T}^j) - \prod_{j=1}^d K_h(x^j - X_t^j(\frac{t}{T}))
| \le C \sum_{k=1}^d | K_h(x^k - X_{t,T}^k) - K_h(x^k -
X_t^k(\frac{t}{T})) |$. Once again exploiting the boundedness of $K$,
we can find a constant $C< \infty$ with $| K_h(x^k - X_{t,T}^k) -
K_h (x^k - X_t^k(\frac{t}{T})) | \le C | K_h(x^k - X_{t,T}^k) -
K_h(x^k - X_t^k(\frac{t}{T})) |^r$ for $r = \min\{\rho,1\}$. Hence,
%
%eC.1 #&#
%
\begin{eqnarray}\label{kernel-product}
&& \Biggl| \prod_{j=1}^d K_h
\bigl(x^j - X_{t,T}^j\bigr) - \prod
_{j=1}^d K_h \biggl(x^j -
X_t^j\biggl( {\frac{t}{T}}\biggr) \biggr) \Biggr|
\nonumber\\[-9pt]\\[-9pt]
&&\qquad\le C \sum_{k=1}^d \biggl| K_h
\bigl(x^k - X_{t,T}^k\bigr) - K_h
\biggl(x^k - X_t^k\biggl( {\frac
{t}{T}}
\biggr) \biggr) \biggr|^r.
\nonumber
\end{eqnarray}
Using (\ref{kernel-product}), we obtain
\begin{eqnarray*}
&& \bigl|Q_1(u,x)\bigr|
\\[-2pt]
&&\qquad\le\frac{C}{Th^{d+1}} \sum_{t=1}^T
K_h \biggl(u - \frac
{t}{T} \biggr) \\[-2pt]
&&\qquad\quad\hspace*{47pt}{}\times\ex\Biggl[ \sum
_{k=1}^d \biggl| K_h\bigl(x^k -
X_{t,T}^k\bigr) - K_h \biggl(x^k -
X_t^k\biggl( {\frac{t}{T}}\biggr) \biggr)
\biggr|^r
\\[-2pt]
&&\qquad\quad\hspace*{70.5pt}{} \times\prod_{j=1}^d
\bar{K}_h\bigl(x^j - X_{t,T}^j\bigr)
\biggl| m \biggl(\frac{t}{T},X_{t,T} \biggr) - m(u,x) \biggr| \Biggr]
\end{eqnarray*}
with $r = \min\{\rho,1\}$. The term $\prod_{j=1}^d \bar
{K}_h(x^j - X_{t,T}^j) | m(\frac{t}{T},X_{t,T}) - m(u,x) |$ in the
above expression can be bounded by $Ch$. Since $K$ is Lipschitz,
$|X_{t,T}^k - X_t^k(\frac{t}{T})| \le\frac{C}{T} U_{t,T}(\frac
{t}{T})$ and the variables $U_{t,T}(\frac{t}{T})$ have finite $r$th
moment, we can infer that
\begin{eqnarray*}
&& \bigl|Q_1(u,x)\bigr|
\\
&&\qquad\le\frac{C}{Th^d} \sum_{t=1}^T
K_h \biggl(u - \frac
{t}{T} \biggr) \ex\Biggl[ \sum
_{k=1}^d \biggl| K_h\bigl(x^k -
X_{t,T}^k\bigr) - K_h \biggl(x^k -
X_t^k\biggl( {\frac{t}{T}}\biggr) \biggr)
\biggr|^r \Biggr]
\\
&&\qquad\le\frac{C}{Th^d} \sum_{t=1}^T
K_h \biggl(u - \frac
{t}{T} \biggr) \ex\Biggl[ \sum
_{k=1}^d \biggl| \frac{1}{Th} U_{t,T}
\biggl( {\frac{t}{T}}\biggr) \biggr|^r \Biggr] \\
&&\qquad\le{\frac{C}{T^r h^{d-1+r}}}
\end{eqnarray*}
uniformly in $u$ and $x$. Using similar arguments, we can further show
that $\sup_{u,x}|Q_2(u,x)| \le\frac{C}{T^r h^d}$ and $\sup
_{u,x}|Q_3(u,x)| \le\frac{C}{T^r h^{d-1+r}}$. Finally, applying
Lemmas~\ref{lemmaB1} and~\ref{lemmaB4} and exploiting the smoothness
conditions on $m$ and~$f$, we obtain that uniformly in $u$ and $x$,
\[
Q_4(u,x) = h^2 \frac{\kappa_2}{2} \sum
_{i=0}^d \bigl( 2 \,\partial_i m(u,x)\,
\partial_i f(u,x) + \partial_{ii}^2 m(u,x)
f(u,x) \bigr) + o\bigl(h^2\bigr).
\]
Combining the results on $Q_1(u,x),\ldots,Q_4(u,x)$ yields~\ref{conv-c}.
\end{pf*}

\subsection*{Proof of Theorem \protect\ref{theo-normality}}

The result can be shown by using the techniques from Theorem \ref
{theo-nw} together with a blocking argument. More details are given in
the supplement~\cite{Vogt2012}.

%sD #&#
\section{}\label{appC}

In this Appendix, we prove the results concerning the smooth
backfitting estimates of Section~\ref{section-add}. Throughout the
Appendix, conditions~\ref{add1} and~\ref{add2} are assumed
to be satisfied.

\subsection*{Auxiliary results}

Before we come to the proof of Theorems~\ref{theorem-sbf-1} and \ref
{theorem-sbf-2}, we provide results on uniform convergence rates for
the kernel smoothers that are used as pilot estimates in the smooth
backfitting procedure. We start with an auxiliary lemma which is needed
to derive the various rates.
%
%leD.1 #&#
%
\begin{lemma}\label{lemmaC-obs}
Define $T_0= \ex[\Obs]$. Then uniformly for $u \in I_h$,
%
%eD.1 #&#
%
\begin{equation}
\label{obs1} \frac{T_0}{T} = \pr\bigl(X_0(u)
\in[0,1]^d\bigr) + O\bigl(T^{-{\rho}/({1+\rho
})}\bigr) + o(h)
\end{equation}
with $\rho$ defined in assumption~\ref{C1} and
%
%eD.2 #&#
%
\begin{equation}
\label{obs2} \frac{\Obs- T_0}{T_0} = O_p \biggl(\sqrt{
\frac{\log T}{T h}} \biggr).
\end{equation}
\end{lemma}
\begin{pf}
The proof can be found in the supplement
\cite{Vogt2012}.
\end{pf}

We now examine the convergence behavior of the pilot estimates of the
backfitting procedure. We first consider the density estimates $\hat
{p}_j$ and $\hat{p}_{j,k}$.
%
%leD.2 #&#
%
\begin{lemma}\label{lemma-add-kde}
Define $v_{T,2} = \sqrt{\log T / T h^2}$, $v_{T,3} = \sqrt{\log T / T
h^3}$ and $b_{T,r} = T^{-r} h^{-(d+r)}$ with $r = \min\{ \rho,1 \}$.
Moreover, let $\kappa_0(w) = \int K_h(w,v)\,dv$. Then
\begin{eqnarray*}
\sup_{u,x^j \in I_h} \bigl| \hat{p}_j\bigl(u,x^j\bigr) -
p_j\bigl(u,x^j\bigr) \bigr| & = & O_p(v_{T,2})
+ O(b_{T,r}) + o(h),
\\
\sup_{u \in I_h, x^j \in[0,1]} \bigl| \hat{p}_j\bigl(u,x^j\bigr) -
\kappa_0\bigl(x^j\bigr) p_j
\bigl(u,x^j\bigr) \bigr| & = & O_p(v_{T,2}) +
O(b_{T,r}) + O(h),
\\
\sup_{u,x^j,x^k \in I_h} \bigl| \hat{p}_{j,k}\bigl(u,x^j,x^k
\bigr) - p_{j,k}\bigl(u,x^j,x^k\bigr) \bigr| & = &
O_p(v_{T,3}) + O(b_{T,r}) + o(h),\\[-35pt]
\end{eqnarray*}

\begin{eqnarray*}
&&\qquad\quad\hspace*{5pt}\mathop{\sup_{u \in I_h, }}_{ x^j,x^k \in[0,1]} \bigl| \hat{p}_{j,k}
\bigl(u,x^j,x^k\bigr) - \kappa_0
\bigl(x^j\bigr) \kappa_0\bigl(x^k\bigr)
p_{j,k}\bigl(u,x^j,x^k\bigr) \bigr|
\\
&&\hspace*{140.5pt}\qquad\qquad = O_p(v_{T,3}) + O(b_{T,r}) + O(h).
\end{eqnarray*}
\end{lemma}
\begin{pf}
We only consider the term $\hat{p}_j$, the proof for
$\hat{p}_{j,k}$ being analogous. Defining
$\check{p}_j(u,x^j) = (T_0)^{-1} \sum_{t=1}^T I(X_{t,T} \in
[0,1]^d) K_h(u,\frac{t}{T}) K_h(x^j,X_{t,T}^j)$
with $T_0= \ex[\Obs]$, we obtain that
\begin{eqnarray*}
\hat{p}_j\bigl(u,x^j\bigr) & = & \biggl[ 1 +
\frac{\Obs- T_0}{T_0} \biggr]^{-1} \check{p}_j
\bigl(u,x^j\bigr)
\\
& = & \biggl[ 1 - \frac{\Obs- T_0}{T_0} + O_p \biggl(\frac{\Obs-
T_0}{T_0}
\biggr)^2 \biggr] \check{p}_j\bigl(u,x^j
\bigr).
\end{eqnarray*}
By (\ref{obs2}) from Lemma~\ref{lemmaC-obs}, this implies that
\[
\hat{p}_j\bigl(u,x^j\bigr) = \check{p}_j\bigl(u,x^j\bigr) + O_p(\sqrt{\log T / T h})
\]
uniformly for $u \in I_h$ and $x^j \in[0,1]$. Applying the proving
strategy of Theorem~\ref{theo-nw} to $\check{p}_j(u,x^j)$ completes
the proof.
\end{pf}

We next examine the Nadaraya--Watson smoother $\hat{m}_j$. To this
purpose, we decompose it into a variance part $\hat{m}_j^A$ and a bias
part $\hat{m}_j^B$. The decomposition is given by $\hat{m}_j(u,x^j) =
\hat{m}_j^A(u,x^j) + \hat{m}_j^B(u,x^j)$ with
\begin{eqnarray*}
\hat{m}_j^A\bigl(u,x^j\bigr) &=&
\frac{1}{\Obs} \sum_{t=1}^T I
\bigl(X_{t,T} \in[0,1]^d\bigr)
K_h \biggl(u,
\frac{t}{T} \biggr) K_h\bigl(x^j,X_{t,T}^j
\bigr) \varepsilon_{t,T} / \hat{p}_j\bigl(u,x^j
\bigr),
\\
\hat{m}_j^B\bigl(u,x^j\bigr) & = &
\frac{1}{\Obs} \sum_{t=1}^T I
\bigl(X_{t,T} \in[0,1]^d\bigr) K_h \biggl(u,
\frac{t}{T} \biggr) K_h\bigl(x^j,X_{t,T}^j
\bigr)
\\
&&\hspace*{44.3pt}{} \times\Biggl( m_0 \biggl(\frac{t}{T} \biggr) + \sum
_{k=1}^d m_k \biggl(
\frac{t}{T},X_{t,T}^k \biggr) \Biggr) \Big/
\hat{p}_j\bigl(u,x^j\bigr).
\end{eqnarray*}
The next two lemmas characterize the asymptotic behavior of $\hat
{m}_j^A$ and $\hat{m}_j^B$.\vadjust{\goodbreak}

%leD.3 #&#
%
\begin{lemma}\label{lemma-add-nw-var}
It holds that
%
%eD.3 #&#
%
\begin{equation}
\label{sbf-nw-var} \sup_{u, x^j \in[0,1]} \bigl| \hat{m}_j^A
\bigl(u,x^j\bigr) \bigr| = O_p \biggl( \sqrt{
\frac{\log T}{T h^2}} \biggr).
\end{equation}
\end{lemma}
\begin{pf}
Replacing the occurrences\vspace*{1pt} of $\Obs$ in $\hat
{m}_j^A$ by $T_0= \ex[\Obs]$ and then applying Theorem \ref
{theo-stochastic-part} gives the result.
\end{pf}
%
%leD.4 #&#
%
\begin{lemma}\label{lemma-add-nw-bias}
It holds that
%
%eD.4 #&#
%eD.5 #&#
%
\begin{eqnarray}
\label{bias-expansion1}
\sup_{u, x^j \in I_h} \bigl| \hat{m}_j^B\bigl(u,x^j
\bigr) - \hat{\mu}_{T,j}\bigl(u,x^j\bigr) \bigr| & = &
o_p\bigl(h^2\bigr),
\\
\label{bias-expansion2}
\sup_{u \in I_h, x^j \in I_h^c} \bigl| \hat{m}_j^B\bigl(u,x^j
\bigr) - \hat{\mu}_{T,j}\bigl(u,x^j\bigr) \bigr| & = &
O_p\bigl(h^2\bigr)
\end{eqnarray}
with $I_h^c = [0,1] \setminus I_h$ and
\begin{eqnarray*}
\hat{\mu}_{T,j}\bigl(u,x^j\bigr) & = & \alpha_{T,0}(u)
+ \alpha_{T,j}\bigl(u,x^j\bigr)\\
&&{} + \sum
_{k \ne j} \int\alpha_{T,k}\bigl(u,x^k\bigr)
\frac{\hat
{p}_{j,k}(u,x^j,x^k)}{\hat{p}_j(u,x^j)} \,dx^k
\\
&&{} + h^2 \int\beta(u,x) \frac{p(u,x)}{p_j(u,x^j)} \,dx^{-j}.
\end{eqnarray*}
Here,
\begin{eqnarray*}
\alpha_{T,0}(u) & = & m_0(u) + h \kappa_1(u)
\,\partial_u m_0(u) + \frac{h^2}{2}
\kappa_2(u) \,\partial_{uu}^2
m_0(u),
\\
\alpha_{T,k}\bigl(u,x^k\bigr) & = & m_k
\bigl(u,x^k\bigr) + h \biggl[ \kappa_1(u)
\,\partial_u m_k\bigl(u,x^k\bigr) +
\frac
{\kappa_0(u) \kappa_1(x^k)}{\kappa_0(x^k)} \,\partial_{x^k} m_k\bigl(u,x^k
\bigr) \biggr],
\\
\beta(u,x) & = & \kappa_2 \,\partial_u m_0(u)
\,\partial_u \log p(u,x)
\\
&&{} + \sum_{k=1}^d \biggl\{
\kappa_2 \,\partial_u m_k\bigl(u,x^k
\bigr) \,\partial_u \log p(u,x) + \frac{\kappa_2}{2}
\,\partial_{uu}^2 m_k\bigl(u,x^k
\bigr)
\\
&&\hspace*{31pt}{} + \kappa_2 \,\partial_{x^k} m_k
\bigl(u,x^k\bigr) \,\partial_{x^k} \log p(u,x) +
\frac{\kappa_2}{2} \,\partial_{x^k x^k}^2 m_k
\bigl(u,x^k\bigr) \biggr\},
\end{eqnarray*}
where the symbol $\partial_z g$ denotes the partial derivative of the
function $g$ with respect to $z$ and $\kappa_2 = \int w^2 K(w) \,dw$ as
well as $\kappa_l(v) = \int w^l K_h(v,w)\,dw$ for $l=0,1,2$.
\end{lemma}
\begin{pf}
As the proof is rather lengthy and involved, we only
sketch its idea. A detailed version can be found in the supplement
\cite{Vogt2012}. To provide the stochastic expansion of $\hat
{m}_j^B(u,x^j)$ in (\ref{bias-expansion1}) and (\ref
{bias-expansion2}), we follow the proving strategy of Theorem 4 in
Mammen et al.~\cite{Mammen1999}. Adapting this strategy is, however,
not completely straightforward. The complication mainly results from
the fact that we cannot work with the variables $X_{t,T}$ directly but
have to replace them by the approximations $X_t(\frac{t}{T})$. To cope
with the resulting difficulties, we exploit (\ref{obs1}) and (\ref
{obs2}) of Lem\-ma~\ref{lemmaC-obs} and use arguments similar to those
for Theorem~\ref{theo-nw}.
\end{pf}

We finally state a result on the convergence behavior of the term
$\tilde{m}_0(u)$.
%
%leD.5 #&#
%
\begin{lemma}\label{lemma-add-m0}
It holds that
%
%eD.6 #&#
%
\begin{equation}
\label{sbf-m0} \sup_{u \in I_h} \bigl| \tilde{m}_0(u) -
m_0(u) \bigr| = O_p \biggl( \sqrt{\frac{\log T}{T h}} +
h^2 \biggr).
\end{equation}
\end{lemma}
\begin{pf}
The claim can be shown by replacing $\Obs$ with
$T_0= \ex[\Obs]$ in the expression for $\tilde{m}_0(u)$ and then
using arguments from Theorem~\ref{theo-nw}.~%
\end{pf}

\subsection*{Proof of Theorems \protect\ref{theorem-sbf-1} and
\protect\ref{theorem-sbf-2}}

Using the auxiliary results from the previous subsection, it is not
difficult to show that the high-level conditions (A1)--(A6), (A8) and
(A9) of Mammen et al.~\cite{Mammen1999} are satisfied. We can thus
apply their Theorems~\mbox{1--3}, which imply the statements of Theorems \ref
{theorem-sbf-1} and~\ref{theorem-sbf-2}. Note that the high-level
conditions are satisfied uniformly for $u \in I_h$ rather than only
pointwise. Inspecting the proofs of Theorems 1--3 in \cite
{Mammen1999}, this allows us to infer that the convergence rates in
(\ref{rates1}) hold uniformly over $u \in I_h$ rather than only
pointwise. A list of the high-level conditions together with the
details of the proof can be found in the supplement~\cite{Vogt2012}.
\end{appendix}

\section*{Acknowledgments}

%I would like to thank Enno Mammen, Oliver Linton and Suhasini Subba Rao
%for numerous helpful suggestions and comments.

I am grateful to Enno Mammen, Oliver Linton and Suhasini Subba Rao for
numerous helpful suggestions and comments. Moreover, I would like to
thank an Associate Editor and three anonymous referees for their
constructive comments which helped a lot to improve the paper.

\begin{supplement}%[id=suppA]
\stitle{Additional technical details}
\slink[doi]{10.1214/12-AOS1043SUPP} %[doi,text={...}] - jei reikia
%suskaldyti doi
\sdatatype{.pdf}
\sfilename{aos1043\_supp.pdf}
\sdescription{The proofs and technical details that are omitted in the
\hyperref[app]{Appendices} are provided in the supplement that accompanies the paper.}
\end{supplement}

% imsref loaded by lrinkeviciute, 2012-11-14 11:02:44
% imsref loaded by lrinkeviciute, 2012-11-14 11:07:23
%

\printaddresses

\end{document}